\documentclass[10.5pt]{article}
\usepackage{amsmath}
\usepackage{amsfonts}
\usepackage{graphicx}
\usepackage{amsmath}
\usepackage{amssymb}
\usepackage{latexsym}
\usepackage{amsmath,amsfonts,amssymb,amsthm,euscript,makeidx,color,mathrsfs}
\usepackage{hyperref}%%%公式产生超链接
\usepackage{cite} %%%参考文献连续引用时，用破折号连接

\hypersetup{hypertex=true,
colorlinks=true,
linkcolor=blue,
anchorcolor=blue,
citecolor=blue}

\def\sqr#1#2{{\vcenter{\vbox{\hrule height.#2pt
              \hbox{\vrule width.#2pt height#1pt \kern#1pt \vrule width.#2pt}
              \hrule height.#2pt}}}}
\def\signed #1{{\unskip\nobreak\hfil\penalty50
              \hskip2em\hbox{}\nobreak\hfil#1
              \parfillskip=0pt \finalhyphendemerits=0 \par}}
\def\endpf{\signed {$\sqr69$}}
\def\3n{\negthinspace \negthinspace \negthinspace }
\def\2n{\negthinspace \negthinspace }
\def\1n{\negthinspace }

\def\dbE{\mathbb{E}}
\def\dbF{\mathbb{F}}

\def\dbR{\mathbb{R}}
\def\dbS{\mathbb{S}}

\def\sE{\mathscr{E}}

\def\={\buildrel \triangle \over =}

\def\ds{\displaystyle}

\def\ns{\noalign{\ss}}

\def\a{\alpha}
\def\b{\beta}
\def\g{\gamma}
\def\d{\delta}
\def\e{\varepsilon}
\def\z{\zeta}
\def\k{\kappa}
\def\l{\lambda}
\def\m{\mu}

\def\si{\sigma}
\def\t{\tau}
\def\f{\varphi}
\def\th{\theta}

\def\i{\infty}
\def\ti{\tilde}
\def\D{\Delta}
\def\Th{\Theta}
\def\L{\Lambda}

\def\Om{\Omega}

\def\cB{{\cal B}}

\def\cF{{\cal F}}
\def\cG{{\cal G}}
\def\cH{{\cal H}}

\def\cK{{\cal K}}

\def\cP{{\cal P}}
\def\cQ{{\cal Q}}

\def\cX{{\cal X}}

\def\no{\noindent}

\def\ss{\smallskip}
\def\ms{\medskip}

\def\q{\quad}
\def\qq{\qquad}

\def\ra{\rightarrow}

\def\lan{\mathop{\langle}}
\def\ran{\mathop{\rangle}}

\def\pa{\partial}

\def\cd{\cdot}
\def\cds{\cdots}

\def\ae{\hbox{\rm a.e.{ }}}
\def\as{\hbox{\rm a.s.{ }}}

\def\({\Big (}
\def\){\Big )}
\def\[{\Big[}
\def\]{\Big]}

\def\({\Big (}
\def\){\Big )}
\def\[{\Big[}
\def\]{\Big]}

\def\bde{\begin{definition}}
\def\ede{\end{definition}}
\def\be{\begin{equation}}
\def\bel{\begin{equation}\label}
\def\ee{\end{equation}}
\def\bt{\begin{theorem}}
\def\et{\end{theorem}}
\def\bc{\begin{corollary}}
\def\ec{\end{corollary}}
\def\bl{\begin{lemma}}
\def\el{\end{lemma}}
\def\bp{\begin{proposition}}
\def\ep{\end{proposition}}
\def\bas{\begin{assumption}}
\def\eas{\end{assumption}}
\def\br{\begin{remark}}
\def\er{\end{remark}}
\def\ba{\begin{array}}
\def\ea{\end{array}}
\def\bpf{\begin{proof}}
\def\epf{\end{proof}}
\def\ed{\end{document}}

\def\square#1{\vbox{\hrule\hbox{\vrule height#1%
     \kern#1\vrule}\hrule}}
\def\rectangle#1#2{\vbox{\hrule\hbox{\vrule height#1%
     \kern#2\vrule}\hrule}}

\font\tenbb=msbm10 \font\sevenbb=msbm7 \font\fivebb=msbm5

\newfam\bbfam
\scriptscriptfont\bbfam=\fivebb \textfont\bbfam=\tenbb
\scriptfont\bbfam=\sevenbb

\newtheorem{lemma}{Lemma}[section]
\newtheorem{remark}{Remark}[section]

\newtheorem{theorem}{Theorem}[section]
\newtheorem{corollary}{Corollary}[section]

\newtheorem{definition}{Definition}[section]
\newtheorem{proposition}{Proposition}[section]
\newtheorem{assumption}{Assumption}[section]

\makeatletter
   
   \@addtoreset{equation}{section}
\makeatother

\def\ges{\geqslant}
\def\les{\leqslant}
\def\rf{\eqref}
\def\blan{\big\langle}
\def\bran{\big\rangle}
\def\Blan{\Big\langle}
\def\Bran{\Big\rangle}

\begin{document}

\title{A general maximum principle for optimal control of stochastic differential delay systems\thanks{
{\bf Funding}: {This work was supported by the National Key $R\&D$ Program of China (2018YFA0703800, 2022YFA1006104), the National Natural Science Foundation of China (T2293770, 11971266, 11831010, 11971332, 11931011), the Shandong Provincial Natural Science Foundations (ZR2022JQ01, ZR2020ZD24, ZR2019ZD42), and the Science Development Project of Sichuan University (2020SCUNL201).}}}

\author{Weijun Meng\thanks{ Key Laboratory of Systems and Control, Institute of Systems Science, Academy of Mathematics and Systems Science, Chinese Academy of Sciences,  Beijing 100190, China. ({mengwj@mail.sdu.edu.cn}).}
\and Jingtao Shi\thanks{ School of Mathematics, Shandong University, Jinan 250100, China. ({shijingtao@sdu.edu.cn}). }
 \and Tianxiao Wang\thanks{ Corresponding author. School of Mathematics, Sichuan University, Chengdu 610065, China. ({wtxiao2014@scu.edu.cn}).}
\and Ji-Feng Zhang\thanks{Key Laboratory of Systems and Control, Institute of Systems Science,
Academy of Mathematics and Systems Science, Chinese Academy of Sciences, Beijing 100190, China,
and School of Mathematics Sciences, University of Chinese Academy of Sciences, Beijing 100149,
China. ({jif@iss.ac.cn}).}}

\maketitle

\begin{abstract}
In this paper, we solve an open problem and obtain a general maximum principle for a stochastic optimal control problem where the control domain is an arbitrary non-empty set and all the coefficients (especially the diffusion term and the terminal cost) contain the control and state delay. In order to overcome the difficulty of dealing with the cross term of state and its delay in the variational inequality, we propose a new method: transform a delayed variational equation into a Volterra integral equation without delay, and introduce novel first-order, second-order adjoint equations via the backward stochastic Volterra integral equation theory. Finally we express these two kinds of adjoint equations in more compact anticipated backward stochastic differential equation types for several special yet typical control systems.
\end{abstract}

\bf Keywords. \rm stochastic differential delay systems, general maximum principle, backward stochastic Volterra integral equations, second-order adjoint equations, non-convex control domain

\bf AMS Mathematics subject classification. \rm 93E20, 60H20, 34K50

\section{Introduction}

The study of optimal control problem has been a hot topic for decades, and maximum principle has been one of the main approaches to address the control problems. In 1965, Kushner (see \cite{Kushner-1965}) firstly studied the maximum principle for the stochastic optimal control problem, where the diffusion term does not contain state and control. Since then, extensive literatures have emerged to study the stochastic optimal control problems. However, either the control domain must be convex, or the diffusion term does not contain the control. In 1990, Peng (see \cite{Peng-1990}) completely solved the stochastic optimal control problem and obtained the general maximum principle, by means of \emph{backward stochastic differential equations} (BSDEs) as adjoint equations. On the other hand, in the real world, the memory affect always exists. The increment of the control system  not only depends on the current state, but also depends on the past state.  Also when the controller decides to exert control, it take some time to exercise the action. Therefore, it has profound theory importance and extensive application value to study the control problems for systems with both state delay and control delay. Usually \emph{stochastic differential delay equations} (SDDEs) are used to describe these delayed control systems. More details about SDDEs can be referred to \cite{Kushner08,Mao97,Mohammed-1984,Mohammed98}.

\ss

Given a time duration $[0,T]$, for a non-empty set $U\subset \dbR^m$, not necessarily convex, a constant time delay parameter $\d\in(0,T)$ and a constant $\l\in\dbR$, in this paper we consider the system of the following form:
\vskip-7mm
\bel{state-}\left\{\ba{ll}
\ds d x(t)=b\(t,x(t),x(t-\d),\int_{-\d}^0e^{\l\th}x(t+\th)d\th,u(t),u(t-\d)\)dt\\
\ns\ds\q \qq+\si\(t,x(t),x(t\1n-\1n\d),\2n\int_{-\d}^0\3ne^{\l\th}x(t+\th)d\th,\1nu(t), \1nu(t\1n-\1n\d)\)dW(t),t\in[0,T],\\
\ns\ds x(t)=\xi(t),\ u(t)=\eta(t),\ t\in[-\d,0],
\ea\right.\ee
\vskip-1mm
\no where $x(\cd)\in\dbR^n$ is state and $u(\cd)\in U$ is control. Suppose that $(\Omega,\mathcal{F},\dbF,\mathbb{P})$ is a complete filtered probability space and the filtration $\dbF=\{\mathcal{F}_t\}_{t\geq 0}$ is generated by a $d$-dimensional standard Brownian motion $\{W(t)\}_{t\geq0}$. $b,\si$ are given random coefficients with proper dimensions. Deterministic continuous function $\xi(\cd)$ and square integrable function $\eta(\cd)$ are the initial trajectories of the state and the control, respectively. We associate \rf{state-} with the following cost functional
\vskip-3mm
\bel{cost-}\ba{ll}
\ns\ds J(u(\cd))=\dbE\bigg[\int_0^T l\(t,x(t),x(t-\d),\int_{-\d}^0e^{\l\th}x(t+\th)d\th,u(t),u(t-\d)\)dt \\
\ns\ds \qq\qq\qq+h\(x(T),x(T-\d),\int_{-\d}^0e^{\l\th}X(T+\th)d\th\)\bigg],
\ea\ee
\vskip-1mm
\no where $l,h$ are given random coefficients with proper dimensions.
%$b:[0,T]\times\Om\times\dbR^n\times\dbR^n\times\dbR^n\times U\times U \rightarrow\dbR^n$, $\sigma:[0,T]\times\Om\times\dbR^n\times\dbR^n\times\dbR^n\times U\times U \rightarrow\dbR^{n\times d}$, $l:[0,T]\times\Om\times\dbR^n\times\dbR^n\times\dbR^n\times U\times U \rightarrow\dbR$ and $h:\Om\times\dbR^n\times\dbR^n\times\dbR^n\rightarrow\dbR$ are given mappings.
%
Define the admissible control set as follows:
\vskip-7mm
$$\ba{ll}
\ns\ds \mathcal{U}_{ad}:=\Big\{u(\cdot):[-\d,T]\ra\dbR^m\big|u(\cdot)\mbox{ is a }U\mbox{-valued, square-integrable},\ \dbF\mbox{-adapted}\\
\ns\ds  \qq\qq\qq\qq\qq\qq\q \mbox{process and }u(t)=\eta(t),\ t\in[-\d,0]\Big\}.
\ea$$
\vskip-3mm

We state the optimal control problem as follows:

\textbf{Problem (P).} Our object is to find a control $u^*(\cdot)$ over $\mathcal{U}_{ad}$ such that (\ref{state-}) is satisfied and (\ref{cost-}) is minimized, i.e.,
$$\ba{ll}
J(u^*(\cd))=\inf\limits_{u(\cd)\in\,\mathcal{U}_{ad}}J(u(\cd)).
\ea$$
\vskip-1mm

Any $u^*(\cdot)\in \mathcal{U}_{ad}$ that achieves the above infimum is called an {\it optimal control} and the corresponding solution $x^*(\cdot)$ is called the {\it optimal trajectory}. $(x^*(\cdot),u^*(\cdot))$ is called an {\it optimal pair}. Optimal control problems of stochastic differential delay systems are widely used in economics, engineering and medicine (see \cite{Chen-Wu-2020,Mao-Sabanis-2013,Shen-Zeng-2014, Xu-Shi-Zhang-2018}), and thus have attracted more and more scholars' attention. Take an optimal consumption problem as an example, at time $t$ let $x(t)$, $u(t)$ be the wealth, the consumption amount, respectively. It is reasonable to suppose that the wealth increment is a combination of the present value $x(t)$ plus some sliding average of previous value $\int_{-\d}^0e^{\l\th}x(t+\th)d\th$ and negative consumption amount $u(t)$. Therefore, the wealth equation satisfied by $x(\cd)$ has the form of \rf{state-}. The consumer always wants to find an optimal consumption strategy $u^*(\cd)$ to maximize his terminal wealth $\dbE[X(T)]$ and consumption satisfaction $\dbE\int_0^T\frac{u^\g(t)}{\g}dt$, where $\g\1n\in\1n(0,\1n1)$, $1-\g$ is the relative risk aversion of the consumer. Thus, the cost functional \rf{cost-} can be chosen as $\dbE\big[\1n-\1nX(T)\1n-\1n\int_0^T\1n\frac{u^\g(t)}{\g}\1ndt\big]$. With different levels of consumption packages for consumers to select, the value set $U$ of the consumption amount $u(t)$ should be limited and not necessarily convex. This typical consumption problem is a case of Problem (P), which motivates us to study the maximum principle for Problem (P).

\ss

So far, there have been extensive literatures to study optimal control problems of stochastic differential delay systems. \O ksendal and Sulem in \cite{Oksendal-Sulem-2000} studied the sufficient maximum principle for the stochastic optimal control problem with convex control domain, and required the solution of certain adjoint equation to be zero due to the lack of It\^o formula to deal with pointwise state delay terms. Chen and Wu in \cite{Chen-Wu-2010} introduced a class of anticipated BSDEs as the adjoint equations and obtained the maximum principle. Although \cite{Chen-Wu-2010} removed the ``zero-solution" condition in \cite{Oksendal-Sulem-2000}, the control domain is still convex. Recently, Meng and Shi in \cite{Meng-Shi-2021} addressed the stochastic optimal control problem, allowed the control domain to be non-convex, and gave the general maximum principle. However, the solution of some second-order adjoint equation must be zero, since at that moment there is  no proper method to eliminate the cross terms of states and their delay terms. More related literatures can be referred to \cite{Chen-Yu-2015,Guatteri-Masiero-2020,Huang-Shi-2012, Ni-Yiu-Zhang-Zhang-2017,Wu-Wang-2015,Yu-2012,Zhang-2021, Zhang-2022,Zhang-Xu-2017}.

\ss

In this paper, we consider the stochastic optimal control problem associated with (\ref{state-}), (\ref{cost-}), and derive the general maximum principle with arbitrary non-empty control domain $U$. Different from all the aforementioned literatures, we study the optimal control problem from a new viewpoint of forward stochastic Volterra integral systems and develop some effective techniques. More precisely, inspired by \cite{Yushi-delay}, we first properly transform the delayed first-order variational equation into a linear forward \emph{stochastic Volterra integral equation} (SVIE) without delay. Then, we combine it with the original first-order variational equation, lift them up, and end up with a higher dimensional linear forward SVIE. Eventually, we adapt the arguments developed by Wang and Yong (see \cite{Wang-Yong-2022}) for optimal control problems of forward stochastic Volterra integral systems into our framework and derive the main results accordingly.

\ss

Forward Volterra integral systems were introduced by Italian mathematician Volterra (see \cite{Volterra-book-1930}). So far there have been extensive literatures about the optimal control problems of forward Volterra integral systems. However, there are very little work to study the optimal control of forward stochastic Volterra integral systems. One possible reason is that until 2002 the theory of Type-I backward SVIEs was established by Lin (see \cite{Lin-2002}). Then, in 2006 Yong (see \cite{Yong-2008}) proposed Type-II backward SVIEs and firstly derived the maximum principle for optimal control problems of forward stochastic Volterra integral systems  with convex control domain. Until recently, Wang and Yong in \cite{Wang-Yong-2022} introduced an auxiliary process and obtained the general maximum principle, where the control domain is allowed to be non-convex.  More references can be referred to \cite{Wang-2020,Wang-Zhang-2017}.

As far as we know, a number of papers transform the delayed control problem into a control problem of Volterra integral systems. For example, in \cite{Huang-Li-Wang-2016}, they used proper variation of constants formula to transform equivalently the delayed quadratic optimal control problem into that of a linear Volterra integral system. Similar ideas also happened in \cite{Lee-You-1989} in infinite dimensional setting.
On the other hand, there are also other methods to transform the delayed system to another system (see  \cite{Curtain-Pritchard-1978,Delfour-Mitter-1972,Duan-2023, Ichikawa-1982,Vinter-Kwong-1981}). Among them, the delayed finite dimensional problem was lifted up to an infinite dimensional problem without delay. A limitation of such method lies in the high regularity assumption (such as continuity and differentiability) for the coefficients when going back to the original problem. Notice that our transformation in the current paper are essentially different from the above. In addition, by our arguments on (\ref{state-}), there is no need to introduce infinite dimensional analysis.

\ss

The innovations and contributions of this paper are as follows:

{\bf(i)} The control system is very general. The control domain is not required to be convex, pointwise and distributed state delays appear not only in the state equation and the running cost, but also in the terminal cost which is new in the existing literatures to our best knowledge, because there exist many difficulties in seeking adjoint equations for variational equations of pointwise state delays, and pointwise control delays can appear in the diffusion term and the running cost. Thus, our model can cover most control systems in the existing literatures, such as \cite{Chen-Wu-2010,Meng-Shi-2021,Oksendal-Sulem-2000,Zhang-Xu-2017}.

{\bf(ii)} A general maximum principle is obtained. It is simple and concise, consisting of two parts: one describes the maximum condition with delay, and the other describes the maximum condition without delay. In contrast with \cite{Meng-Shi-2021}, the strict ``zero-solution" condition imposed on the adjoint equation is successfully removed.

{\bf(iii)} A new method is proposed to treat cross terms. How to deal with the cross terms ``$ x_1(t)^\top[\cds]y_1(t)$" and ``$ y_1(t)^\top[\cds]x_1(t)$" in the variational inequality, is a key yet difficult problem in obtaining the general maximum principle. Inspired by \cite{Wang-Yong-2022}, we solve this hard issue by the theory of forward, backward stochastic Volterra integral systems.

{\bf(iv)} Novel adjoint equations are introduced. The first-order adjoint equations consist  of a simple BSDE and a backward SVIE, while the second-order adjoint equations consist of a simple BSDE and three coupled backward SVIEs. They are used to eliminate the effect of variational processes in the variational inequality, even if the control domain is non-convex and pointwise state delays appear in both the state equation and the terminal cost.

{\bf(v)} The adjoint equations are expressed in more compact forms. The first-order adjoint equation is written as a set of anticipated BSDEs. The second-order adjoint equation reduces to the classical scenario when our delay system reduces to a stochastic differential system.

\ss

The rest of this paper is organized as follows. In Section 2, some basic results are displayed. In Section 3, the delayed variational equations are transformed into Volterra integral equations without delay, and then the adjoint equations are introduced in Section 4. In Section 5, the maximum principle is stated and some careful analysis on the adjoint equations are spread out. Finally, Section 6 gives the concluding remarks.

\section{Preliminaries}

For any $A,B\in\dbR^{m\times d}$, we define by $\lan A,B\ran=Tr[AB^\top]$ the inner product in $\dbR^{m\times d}$ with norm $|\cd|$, and $\dbS^n$ the set of all $n\times n$ symmetric matrices. Let $\dbE_t[\ \cd\ ]\equiv\dbE[\ \cd\ |\mathcal{F}_t]$ be the conditional expectation with respect to $\mathcal{F}_t$, $t\in[0,T]$, and $I$ is the identity matrix of proper dimensions. For $t\in[0,T]$, denote by $L^2_{\cF_t}(\Omega;\dbR^n)$ the Hilbert space consisting of $\dbR^n$-valued $\mathcal{F}_t$-measurable random variable $\xi$ such that $\mathbb{E}|\xi|^2<\infty$,  by $L_{\dbF}^2(0,T;\dbR^n)$ the Hilbert space consisting of $\mathbb{F}$-adapted process $\phi(\cd)$ such that $\mathbb{E}\int_0^T|\phi(t)|^2dt<\infty$, by $L^2_{\dbF}(\Om;C([0,T];\dbR^n))$ the Banach space consisting of $\dbR^n$-valued $\mathbb{F}$-adapted continuous process $\phi(\cd)$ such that $\mathbb{E}\big[\sup\limits_{0\les t\les T}|\phi(t)|^2\big]<\infty$, and by $L^2(0,T;L_\dbF^2(0,T;\dbR^n))$ the space consisting of $\dbR^n$-valued process $\phi(\cd,\cd):[0,T]^2\times\Om\ra\dbR^n$ such that for almost all $t\in[0,T]$, $\phi(t,\cd)\1n\in\1n L^2_\dbF(0,T;\dbR^n)$, $\dbE\int_0^T\int_0^T|\phi(t,s)|^2dsdt<\i$.

\ms

Consider the following SDDE:
\vskip-4mm
\bel{SDDE}\left\{\ba{ll}
\ds d\ti x(t)=b\(t,\ti x(t),\ti x(t-\d),\int_{-\d}^0e^{\l\th}\ti x(t+\th)d\th\)dt\\
\ns\ds \qq\q+\si\(t,\ti x(t),\ti x(t-\d),\int_{-\d}^0e^{\l\th}\ti x(t+\th)d\th\)dW(t),\ \ t\in[0,T],\\
\ns\ds \ti x(t)=\ti \xi(t),\ t\in[-\d,0],
\ea\right.\ee
\vskip-2mm
%
%\no here
%%
%\vskip-8mm
%%
%\bel{yzm-SDDE}\ba{ll}
%%
%\ns\ds \ti y(t):=\ti x(t-\d),\q \ti z(t):=\int_{-\d}^0e^{\l\th}\ti x(t+\th)d\th,
%%
%\ea\ee
%%
%\vskip-3mm
%%
\no $\d>0$ is the constant delay time, $\l\in\dbR$ is a constant, deterministic continuous function $\xi(\cd)$ is the given initial path of the state, and random coefficients $b,\si$ are given mappings satisfying:\\
{\bf(H1)} There exists a constant $L>0$ such that\\
{\color{white}\q} \qq\qq $|\ti b(t,x,y,z)-\ti b(t,x',y',z')|+|\ti\si(t,x,y,z)-\ti\si(t,x',y',z')|$\\
{\color{white}\q} \qq\qq\qq$\les L(|x-x'|+|y-y'|+|z-z'|),\q \forall\ t\in[0,T],\ x,y,z,x',y',z'\in\dbR^n;$\\
{\bf(H2)} $\sup\limits_{0\les t\les T}\big(|\ti b(t,0,0,0)|+|\ti\si(t,0,0,0)|\big)<\i.$

By standard Picard iteration method we can derive the following result.
\bp\label{prop SDDE}
Suppose {\rm (H1)-(H2)} hold. Then, the SDDE \rf{SDDE} admits a unique solution, and there exists a constant $C>0$ such that for $p\ges 2$,
\vskip-7mm
$$\ba{ll}
\ns\ds \2n\dbE\big[\1n\sup\limits_{0\les t\les T}|\ti x(t)|^p\big]\2n\les \1n C\big[\3n\sup\limits_{-\d\les\th\les0}|\ti\xi(\th)|^p \2n+\1n\dbE\big(\2n\int_0^T\3n|\ti b(s,0,0,0)|ds\big)^p \3n+\1n\dbE\big(\2n\int_0^T\3n|\ti\si(s,0,0,0)|^2ds\big)^{\frac p 2}\big].
\ea$$
\ep

\ms

Let $\dbR^+$ be the space of real numbers not less than zero. Consider the following anticipated BSDE:
\vskip-6mm
\bel{ABSDE}\left\{\ba{ll}
\ns\ds \3n -dY(t)\1n=\1ng\big(t,\1nY(t),\1nZ(t),\1nY(t\1n+\1n\d^1(t)),\1nZ(t\1n+\1n\d^2(t))\big)dt \1n-\2nZ(t)dW(t), t\in[0,T],\\
\ns\ds Y(t)=\alpha(t),\ Z(t)=\beta(t),\q t\in[T,T+K].
\ea\right.\ee
\vskip-2mm
\no Here, terminal conditions $\a(\cdot)\in L_{\dbF}^2(\Om;C([T,T+K];\dbR^m))$ and $\b(\cdot)\in L_{\dbF}^2(T,T+K;\dbR^{m\times d})$ are given, $\d^1(\cdot)$ and $\d^2(\cdot)$ are given $\dbR^+$-valued functions defined on $[0,T]$ satisfying:

\vspace{1mm}

{\bf(H3)} (i) There exists a constant $K\ges0$ such that for all $s\in[0,T],\ s+\d^1(s)\les T+K,\ s+\d^2(s)\les T+K$;

(ii) There exists a constant $M\ges0$ such that for all $t\in[0,T]$ and all nonnegative and integrable function $f(\cdot)$,
\vskip-3mm
\begin{equation*}\begin{aligned}
  \int_t^T\3nf(s+\d^1(s))ds\les M\int_t^{T+K}\3nf(s)ds,\qq \int_t^T\3nf(s+\d^2(s))ds\les M\int_t^{T+K}\3nf(s)ds.
\end{aligned}\end{equation*}
\vskip-2mm

We impose the following conditions to the generator of the equation (\ref{ABSDE}):

\vspace{1mm}

\textbf{(H4)} $g(s,\omega,y,z,\a,\b):\Omega\times\dbR^m\times \dbR^{m\times d}\times L^2_{\mathcal{F}_r}(\Om;\dbR^m)\times L_{\mathcal{F}_{r'}}^2(\Om;\dbR^{m\times d})\rightarrow L^2_{\mathcal{F}_s}\2n(\Om;\dbR^m)$ for all $s\in[0,T]$, where $r,\1nr'\in\1n[s,T\1n+\1nK]$, and $\dbE\big[\int_0^T\2n|g(s,0,0,0,0)|^2ds\big]\1n<\1n+\infty$.

\vspace{1mm}

\textbf{(H5)} There exists a constant $C>0$ such that for all $s\in[0,T]$, $y,\ti y\in\dbR^m$, $z,\ti z\in\dbR^{m\times d}$, $\a(\cd),\ti\a(\cd)\in L_\dbF^2\left(s,T+K;\dbR^m\right)$, $\b(\cd),\ti\b(\cd) \in L_\dbF^2\left(s,T+K;\dbR^{m\times d}\right)$, $r,r^\prime\in[s,T+K]$, we have
\vskip-3mm
$$\ba{ll}
\big|g(s,y,z,\a(r),\b(r^\prime)) -g(s,\ti y,\ti z,\ti\a(r), \ti\b(r^\prime))\big|\\
\ns\ds \qq \les C\big(|y-\ti y| +|z-\ti z|+\dbE_s \big[|\a(r)-\ti\a(r)| +|\b(r^\prime) -\ti\b(r^\prime)|\big]\big).
\ea$$
\vskip-1mm

\bp\label{prop ABSDE} {\rm (see \cite{Peng-Yang-2009})}
Let {\rm (H3)-(H5)} hold. Then, for any given $\a(\cdot)\1n\in\1n L_{\dbF}^2(\Om;\1n$ $C([T,T\1n+\1nK];\dbR^m))$ and $\b(\cdot)\1n\in\1n L_\dbF^2(T,T\1n+\1nK;\dbR^{m\times d})$, the equation (\ref{ABSDE}) admits a unique $\mathcal{F}_t$-adapted solution pair $(Y(\cdot),Z(\cdot))\1n\in\1n L_\dbF^2(\Om;\1nC([0,T\1n+\1nK];\dbR^m))\1n\times\1n L_\dbF^2(0,T\1n+\1nK;\dbR^{m\times d})$.
\ep

\ms

Consider the following backward SVIE:
\vskip-7mm
\bel{BSVIE}\ba{ll}
\ns\ds \ti Y(t)=\psi(t)+\2n\int_t^T\1n\ti g\big(t,s,\ti Y(s),\ti Z(t,s),\ti Z(s,t)\big)ds-\2n\int_t^T\2n\ti Z(t,s)dW(s), t\in[0,T],
\ea\ee
\vskip-3mm
\no where $\ti g$ is the given functions satisfying:\\
%
%$\ti g:\D^*[0,T]\times \dbR^m\times \dbR^{m\times d}\times \dbR^{m\times d}\times\Om\ra\dbR^m$, $\psi:[0,T]\times\Om\ra\dbR^m$
%
{\bf (H6)} $\ti g$ is $\cB([0,T]^2\times\dbR^m\times\dbR^{m\times d}\times \dbR^{m\times d})\otimes\mathcal{F}_T-$measurable such that $s\mapsto\ti g(t,s,y,z,\z)$ is progressively measurable for all $(t,y,z,\z)\in[0,T]\times\dbR^m\times\dbR^{m\times d}\times\dbR^{m\times d}$, and
\vskip-3mm
$$\dbE\int_0^T\bigg(\int_t^T|\ti g(t,s,0,0,0)|ds\bigg)^2dt<\i.$$
\vskip-2mm
\no Moreover,
\vskip-5mm
$$\ba{ll}
\big|\ti g(t,s,y,z,\z)-\ti g(t,s,\bar y,\bar z,\bar\z)\big| \les L(t,s)\big(|y-\bar y|+|z-\bar z|+|\z-\bar\z|\big),\\
\ns\ds \qq\qq\qq\q \forall\ 0\les t\les s\les T,\ y,\bar y\in\dbR^m,\ z,\bar z,\z,\bar \z\in\dbR^{m\times d},\ \as
\ea$$
\vskip-2mm
\no where $L$ is a deterministic function such that for some $\e>0$,
\vskip-3mm
$$\sup\limits_{t\in[0,T]}\int_t^TL(t,s)^{2+\e}ds<\i.$$
\vskip-2mm

\bp\label{pro for bsvie}{\rm (see \cite{Yong-2008})}
Let {\rm (H6)} hold. Then, for any $\cB([0,T])\otimes\mathcal{F}_T$-measurable process $\psi(\cd)$ satisfying $\dbE\2n\int_0^T\2n|\psi(t)|^2dt\1n<\2n\i$, the backward SVIE \rf{BSVIE} admits a unique adapted solution $(Y(\cd),Z(\cd,\cd))\in L_{\dbF}^2(0,T;\dbR^m)\times L^2(0,T;L_{\dbF}^2(0,T;\dbR^{m\times d}))$ satisfying
\vskip-3mm
 $$\ti Y(t)=\dbE_s[\ti Y(t)]+\int_s^t\ti Z(t,r)dW(r),\ \ae t\in[s,T].$$
\vskip-2mm
\no Moreover, for any $s\in[0,T]$, the following estimate holds:
\vskip-4mm
$$\ba{ll}
\ns\ds \dbE\bigg[\int_s^T|\ti Y(t)|^2dt+\2n\int_s^T\2n\int_s^T\2n|\ti Z(t,r)|^2drdt\bigg]
\\
\ns\ds \qq\qq\1n\les \1n C\dbE\bigg[\int_s^T\2n|\psi(t)|^2dt\1n+\1n\int_s^T\3n\bigg(\1n\int_t^T\3n| \ti g(t,r,0,0,0)|dr\bigg)^2dt\bigg].
\ea$$
\ep

\section{A novel transformation from SDDE to SVIE }

In this section, we present the variational equations to be studied, then make some interesting transformations to them. Similar transformation also appeared in \cite{Yushi-delay}.

\ms

Denote
\bel{yzm}\ba{ll}
\ns\ds y(t):=x(t-\d),\q z(t):=\int_{-\d}^0e^{\l\th}x(t+\th)d\th,\q \m(t):=u(t-\d).
\ea\ee
\vskip-3mm
\no Then, we can rewrite the state equation \rf{state-} in a more concise form as follows:
\vskip-4mm
\bel{state}\left\{\ba{ll}
\ns\ds d x(t)=b(t,x(t),y(t),z(t),u(t),\m(t))dt\\
\ns\ds\qq \qq+\si(t,x(t),y(t),z(t),u(t),\m(t))dW(t),\ \ t\in[0,T],\\
\ns\ds x(t)=\xi(t),\ u(t)=\eta(t),\ t\in[-\d,0].
\ea\right.\ee
\vskip-1mm
\no And the cost \rf{cost-} becomes
\vskip-4mm
\bel{cost}\ba{ll}
\ns\ds J(u(\cd))=\dbE\bigg[\int_0^T l(t,x(t),y(t),z(t),u(t),\mu(t))dt+h(x(T),y(T),z(T))\bigg].
\ea\ee
\vskip-1mm

Throughout the paper, we impose the following assumptions.

\textbf{(A1)}
(i) The map $(x,y,z)\mapsto$ $b=b(t,x,y,z,u,\m)$, $\sigma=\sigma(t,x,y,z,u,\m)$, $l=l(t,x,y,z,u,\m)$, $h=h(x,y,z)$ are twice continuously differentiable in $(x,y,z)$. They and all their derivatives $f_{\k^i}$, $f_{\k^i\k^\ell}$ are continuous in $(x,y,z,u,\m)$, $i,\ell=1,2,3$. Here $f=b,\si,l,h$ and $\k^1:=x,\k^2:=y,\k^3:=z$.

(ii) Denote $f=b,\si$ and $g=l,h$. For $i,\ell=1,2,3$, $f_{\k^i}$, $f_{\k^i\k^\ell}$, $g_{\k^i\k^\ell}$ are bounded, where $\k^1:=x,\k^2:=y,\k^3:=z$. There exists a constant $C$ such that
\vskip-7mm
$$
 |f(t,0,0,0,u,\m)|+|g(t,0,0,0,u,\m)|+|g_{\k^i}(t,0,0,0,u,\m)| \les C,\ \forall\  u,\m\in U,\ t\ges 0.
$$
\vskip-2mm

Under (A1), the SDDE \rf{state} admits a unique solution by Proposition \ref{SDDE} above or Theorem 2.1 (\cite{Mohammed-1984}, Chapter II), hence the cost functional \rf{cost} is well-defined and Problem (P) is meaningful.

Since the control domain $U$ is an arbitrary non-empty set, not necessarily convex, we then apply the spike variation technique to deal with Problem (P). Let $u^*(\cdot)$ be the optimal control and $x^*(\cdot)$ be the optimal trajectory. Let $0<\e<\d$, for any given $\t\in[0,T)$, define $u^\varepsilon_\t(t)$ for $t\in[0,T]$ as follows:
\vskip-7mm
\begin{eqnarray}\begin{aligned}\label{perturbed control}
u^\varepsilon_\t(t):=
    \begin{cases}
    u^*(t), & t\notin[\tau,\tau+\varepsilon],\\
    v(t),   & t\in[\tau,\tau+\varepsilon],
    \end{cases}
\end{aligned}\end{eqnarray}
\vskip-2mm
\no which is a perturbed admissible control of the form, where $v(\cdot)$ is any admissible control, and $(x^\varepsilon(\cdot),y^\e(\cd),z^\e(\cd))$ is defined similar to \rf{yzm}. For notional simplicity, for $t\in[0,T]$, we denote $u^\e(t)\equiv u^\e_\t(t)$, $\Theta(t)\equiv(x^*(t),y^*(t),z^*(t)$, $u^*(t)$, $\m^*(t))$, and
\vskip-2mm
\bel{simple notation}\left\{\ba{ll}
\ns\ds f_{\k^i}(t):=f_{\k^i}(t,\Th(t)), \ \ f_{\k^i\k^\ell}(t):=f_{\k^i\k^\ell}(t,\Th(t)),\\
\ns\ds \D f(t):=f(t,x^*(t),y^*(t),z^*(t),u^\e(t),\m^\e(t))-f(t,\Th(t)),\\
\ns\ds \D f_{\k^i}(t):=f_{\k^i}(t,x^*(t),y^*(t),z^*(t),u^\e(t),\m^\e(t)) -f_{\k^i}(t,\Th(t)),\\
\ea\right.\ee
\no where $i,\ell=1,2,3$, $f=b,\si$ and $\k^1:=x,\k^2:=y,\k^3:=z$, $\m^\e(t):=u^\e(t-\d)$. Inspired by \cite{Peng-1990}, we introduce the variational equations:
\vskip-5mm
\bel{variational equation-1}\left\{\ba{ll}
\ns\ds d x_1(t)=\Big[b_x(t)x_1(t)+b_{y}(t)y_1(t)+b_z(t)z_1(t)+\D b(t)\Big]dt\\
\ns\ds \qq\q +\1n\sum_{j=1}^d\2n\Big[\si_x^j(t)x_1(t)\1n+\1n\si_{y}^j(t)y_1(t) \1n+\1n\si_z^j(t)z_1(t)\1n+\1n\D \si^j(t)\Big]\1ndW^j(t), t\1n\in\1n[0,T],\\
\ns\ds x_1(t)=0,\q t\in[-\d,0],
\ea\right.\ee
\vskip-4mm
\bel{variational equation-2}\left\{\ba{ll}
\ns\ds d x_2(t)=\Big[b_x(t)x_2(t)+b_y(t)y_2(t)+b_z(t)z_2(t)\\
\ns\ds \qq\qq +\frac 1 2 \big(x_1(t)^\top,y_1(t)^\top,z_1(t)^\top\big)\pa^2b(t)\big(x_1(t)^\top,y_1(t)^\top,z_1(t)^\top\big)^\top\Big]dt\\
\ns\ds \qq\qq +\sum_{j=1}^d\Big[\si_x^j(t)x_2(t)+\si_y^j(t)y_2(t) +\si_z^j(t)z_2(t) \\
 \ns\ds \qq\qq +\frac 1 2 \big(x_1(t)^\top,y_1(t)^\top,z_1(t)^\top\big)\pa^2\si^j(t)\big(x_1(t)^\top,y_1(t)^\top,z_1(t)^\top\big)^\top\\
\ns\ds\qq\qq +\D \si_x^j(t) x_1(t)+\D \si_y^j(t)y_1(t)+\D\si_z^j(t)z_1(t)\Big]dW^j(t),\q t\in[0,T],\\
\ns\ds x_2(t)=0,\q t\in[-\d,0],
\ea\right.\ee
\vskip-1mm
\no where $y_1(\cd),z_1(\cd),y_2(\cd),z_2(\cd)$ are defined similar to \rf{yzm}, and for $f=b,\si^j$, $j=1,2,\cdots,d$,
\vskip-5mm
$$\ba{ll}
\pa^2f(t)\2n:=\2n\Bigg(\1n\ba{ccc}f_{xx}(t) & f_{xy}(t) & f_{xz}(t)\\
f_{yx}(t) & f_{yy}(t) & f_{yz}(t)\\
f_{zx}(t) & f_{zy}(t) & f_{zz}(t)
\ea\1n\Bigg),
\k^i_1(t)^\top f_{\k^i\k^\ell}(t)\k^\ell_1(t)\2n:=\2n
\begin{pmatrix}\k^i_1(t)^\top  f^1_{\k^i\k^\ell}(t)\k^\ell_1(t)\\
\vdots\\
\k^i_1(t)^\top f^n_{\k^i\k^\ell}(t)\k^\ell_1(t)\end{pmatrix}.
\ea$$
\vskip-1mm
 By Proposition \ref{SDDE}, under Assumption (A1) the variational equations (\ref{variational equation-1}) and (\ref{variational equation-2}) admit a unique solution, respectively. In the following, we introduce some estimates whose proofs are similar to Lemma 3.1 and Lemma 3.2 in \cite{Meng-Shi-2021}.
\begin{lemma}\label{estimate for variational variable}
Let Assumption {\rm (A1)} hold. Then, for any $p\ges1$, we have
$$\ba{ll}
\ns\ds \mathbb{E}\big[\sup\limits_{0\les t\les T}|x^\varepsilon(t)-x^*(t)|^{2p}\big]\1n=\1nO(\varepsilon^p),\q
\mathbb{E}\big[\sup\limits_{0\les t\les T}|x_1(t)|^{2p}\big]\1n=\1nO(\varepsilon^p),\\
\ns\ds \mathbb{E}\big[\sup\limits_{0\les t\les T}|x_2(t)|^p\big]\1n=\1nO(\varepsilon^p),\q
\mathbb{E}\big[\sup\limits_{0\les t\les T}|x^\varepsilon(t)-x^*(t)-x_1(t)|^{2p}\big]\1n =\1no(\varepsilon^p),\\
\ns\ds
\mathbb{E}\big[\sup\limits_{0\les t\les T}|x^\varepsilon(t)-x^*(t)-x_1(t)-x_2(t)|^p\big]\1n=\1no(\varepsilon^p).
\ea$$
\end{lemma}

\bl\label{lemma variational inequality}
\no Let Assumption {\rm (A1)} hold. Suppose $(x^*(\cdot),u^*(\cdot))$ is an optimal pair, $x^\varepsilon(\cd)$ is the trajectory corresponding to $u^\varepsilon(\cdot)$ by (\ref{perturbed control}). Then, the following variational inequality holds:
\bel{variational-inequality}\ba{ll}
\ns\ds  J(u^\e(\cd))\1n-\1nJ(u^*(\cd))\1n=\dbE \Big[h_x(T)\big[x_1(T)+x_2(T)\big]+h_y(T)\big[y_1(T)+y_2(T)\big]\\
\ns\ds \q\1n+h_z(T)\big[z_1(T)+z_2(T)\big]+\frac 1 2\big(x_1(T)^\top\1n,y_1(T)^\top\1n,z_1(T)^\top\1n\big)\pa^2h(T) \\
\ns\ds \q \times\big(x_1(T)^\top\1n,y_1(T)^\top\1n,z_1(T)^\top\1n\big)^\top\Big] +\dbE\int_0^T \Big[\D l(t)+l_x(t)[x_1(t)+x_2(t)] \\
\ns\ds \q+l_y(t)\big[y_1(t)+y_2(t)\big]  +l_z(t)\big[z_1(t)+z_2(t)\big]+\frac 1 2\big(x_1(t)^\top,y_1(t)^\top,z_1(t)^\top\big) \\
\ns\ds \q\times\pa^2l(t)\big(x_1(t)^\top,y_1(t)^\top,z_1(t)^\top\big)^\top\Big]dt+o(\e),
\ea\ee
\vskip-1mm
\no where $\Delta l,l_{\k^i},l_{\k^i\k^\ell},h_{\k^i},h_{\k^i\k^\ell}$ are defined similarly as \rf{simple notation} for $i,\ell=1,2,3$.
\el

\ss

Define
\vskip-4mm
$$\ba{ll}
\ns\ds X_1(t):=\Bigg[\begin{array}{ccccc}
x_1(t) \\
y_1(t){\bf 1}_{(\d,\i)}(t) \\
z_1(t)
\end{array}\Bigg],\q
X_2(t):=\Bigg[\begin{array}{ccccc}
x_2(t) \\
y_2(t){\bf 1}_{(\d,\i)}(t) \\
z_2(t)
\end{array}\Bigg],
\ea$$
\vskip-2mm
\no and for $j=1,\cds,d$,
\vskip-7mm
$$
\ba{ll}
\ns\ds A(t,s):=\Bigg[\begin{array}{ccccc}
b_x(s) & b_y(s) & b_z(s) \\
{\bf 1}_{(\d,\infty)}(t-s) b_x(s) & {\bf 1}_{(\d,\infty)}(t-s) b_y(s) & {\bf 1}_{(\d,\infty)}(t-s) b_z(s)\\
I & -e^{-\l\d}I & -\l I
\end{array}\Bigg],\\
\ns\ds C^j(t,s):=\Bigg[\begin{array}{ccccc}
\si_x^j(s) & \si_y^j(s) & \si_z^j(s) \\
{\bf 1}_{(\d,\infty)}(t-s) \si_{x}^j(s) & {\bf 1}_{(\d,\infty)}(t-s) \si_y^j(s) & {\bf 1}_{(\d,\infty)}(t-s) \si_z^j(s) \\
0 & 0 & 0
\end{array}\Bigg],\\
B(t,s):=\Bigg[\begin{array}{ccccc}
\D b(s) \\
{\bf 1}_{(\d,\infty)}(t-s) \D b (s) \\
0
\end{array}\Bigg],\q
D^j(t,s):=\Bigg[\begin{array}{ccccc}
\D \si^j(s) \\
{\bf 1}_{(\d,\infty)}(t-s) \D \si^j (s) \\
0
\end{array}\Bigg],\\
\ns\ds \bar B(t,\1ns)\2n:=\3n\Bigg[\3n\1n\begin{array}{ccccc}
\frac 1 2 X_1(s)^\top\partial^2 b(s)X_1(s) \\
\frac 1 2 {\bf 1}_{(\d,\infty)}\1n(t\1n-\1ns) \1nX_1(s)^\top\1n\partial^2 b(s)\1nX_1(s) \\
0
\end{array}\3n\1n\Bigg],\q
\D\Xi^j\1n(s)\2n:=\3n\[\1n
\D\si_x^j(s),\D\si_y^j(s),\D\si_z^j(s)\1n\],\\
\ns\ds \bar D^j(t,s):=\Bigg[\begin{array}{ccccc}
\frac 1 2 X_1(s)^\top\partial^2 \si^j(s)X_1(s)+\D \Xi^j (s)X_1(s) \\
{\bf 1}_{(\d,\infty)}(t-s) \big[\frac 1 2 X_1(s)^\top\partial^2 \si^j(s)X_1(s)+\D \Xi^j (s)X_1(s)\big] \\
0
\end{array}\Bigg],
\ea$$
\vskip-3mm
\no where $\bar H$ is $\dbR^{3n}$-valued row vector and other terms are similar. Then, by \rf{variational equation-1}-\rf{variational equation-2},
\vskip-7mm
\bel{variational-1-X}\ba{ll}
\ns\ds X_1(t)\2n=\3n\int_0^t\3n\big[A(t,\1ns)X_1(s)\1n+\1nB(t,\1ns)\big]ds \1n+\sum_{j=1}^d\1n\int_0^t \1n \big[C^j(t,\1ns)X_1(s)\1n+\1nD^j(t,\1ns)\big]dW^j\1n(s),
\ea\ee
\vskip-7mm
\bel{variational-2-X}\ba{ll}
\ns\ds X_2(t)\2n=\3n\int_0^t\3n\big[A(t,\1ns)X_2(s)\1n+\1n\bar B(t,\1ns)\big]ds\1n+\sum_{j=1}^d\1n\int_0^t \1n \big[C^j(t,\1ns)X_2(s)\1n+\1n\bar D^j(t,\1ns)\big]dW^j\1n(s).
\ea\ee
\vskip-2mm

\no Therefore, the above variational inequality (\ref{variational-inequality}) can be written as
\vskip-6mm
\bel{variational-inequality-X}\ba{ll}
\ns\ds J(u^\e(\cd))-J(u^*(\cd))=  \dbE\int_0^T  \Big[\bar L(t)[X_1(t)+X_2(t)] +\frac 1 2  X_1(t)^{\top} L(t) X_1(t)\\
\ns\ds  \qq\q  +\D l(t)\Big]dt+\dbE\[\bar H [X_1(T)+X_2(T)]+\frac 1 2 X_1(T)^{\top} H X_1(T)\]+o(\e).
\ea\ee
\no Here $X_1(\cd)$ and $X_2(\cd)$ satisfy linear SVIEs in (\ref{variational-1-X}) and (\ref{variational-2-X}), respectively, and
$$\ba{ll}
\ns\ds \bar H=\[\begin{array}{ccccc}
h_x(T) & h_y(T) & h_z(T)
\end{array}\],\q \bar L(t)=\[\begin{array}{ccccc}
l_x(t)& l_y(t) & l_z(t)
\end{array}\],\\
\ns\ds H= \Bigg[\begin{array}{ccccc}
h_{xx}(T) & h_{xy}(T) & h_{xz}(T)\\
h_{yx}(T) & h_{yy}(T) & h_{yz}(T)\\
h_{zx}(T) & h_{zy}(T) & h_{zz}(T)
\end{array}\Bigg],\q L(t)=\Bigg[\begin{array}{ccccc}
l_{xx}(t) & l_{xy}(t) & l_{xz}(t)\\
l_{yx}(t) & l_{yy}(t) & l_{yz}(t)\\
l_{zx}(t) & l_{zy}(t) & l_{zz}(t)
\end{array}\Bigg].
\ea$$
\vskip-1mm
\no Under above preparation, we can borrow some useful ideas from \cite{Wang-Yong-2022} where the maximum principle of optimal control problems described by SVIEs was completely solved.

\section{Adjoint equations}
In this section we introduce some adjoint equations to be dual with the variational equations \rf{variational equation-1}-\rf{variational equation-2}.

\subsection{First-order adjoint equations}

We treat the terms about $X_1(\cd)+X_2(\cd)$ in (\ref{variational-inequality-X}).
From \cite{Yong-2008}, we introduce the first-order adjoint equation as follows:
\vskip-3mm
\bel{first-order adjoint equation}\left\{\ba{ll}
\ds (a)\q\eta(t)=\bar H^\top-\sum_{j=1}^d\int_t^T\zeta^j(s) dW^j(s),\q t\in[0,T],\\
\ns\ds (b)\q Y(t)=\1n\bar L(t)^\top\2n+\1nA(T,\1nt)^\top\1n\bar H^\top\2n+\1n\sum_{j=1}^d\1nC^j(T,\1nt)^\top\1n\zeta^j(t) \1n+\2n\int_t^T\3n\Big[\1nA(s,t)^\top\1n Y(s)\\
\ns\ds\qq\qq\q +\1n\sum_{j=1}^d\1nC^j(s,t)^\top\1n Z^j(s,\1nt)\Big]ds \1n-\2n\sum_{j=1}^d\1n\int_t^T\3nZ^j(t,s)dW^j(s), t\1n\in\1n[0,\1nT], \\
\ns\ds (c)\q Y(t)=\dbE Y(t)+\sum_{j=1}^d\int_0^t Z^j(t,s)dW^j(s),\q t\in[0,T].
\ea\right.\ee
\vskip-1mm
\rf{first-order adjoint equation} (a) is a BSDE which admits a unique solution by Theorem 4.1 in \cite{Pardoux-Peng-1990}. On the other hand, \rf{first-order adjoint equation} (b) is a linear backward SVIE, and by Proposition \ref{pro for bsvie}, it admits a unique solution that satisfies
\rf{first-order adjoint equation} (c) under Assumption {\rm (A1)}. Notice that
\vskip-3mm
$$\ba{ll}
\ns\ds X_1(t)\1n+\1nX_2(t)\1n=\1n\f(t)\1n+\1n\int_0^t A(t,s)\big[X_1(s)\1n+\1nX_2(s)\big] ds \1n\\
\ns\ds \qq\qq\qq+\sum_{j=1}^d\int_0^t C^j(t,s)\big[X_1(s)\1n+\1nX_2(s)\big] dW^j(s),
\ea$$
\vskip-1mm
\no where
\vskip-4mm
$$\ba{ll}
\ns\ds \f(t):= \int_0^t\big[\bar B(t,s)+B(t,s)\big]ds+\sum_{j=1}^d\int_0^t \big[\bar D^j(t,s)+D^j(t,s)\big]dW^j(s).
\ea
$$
\vskip-1mm
\no By the dual principle (\cite{Yong-2008}, Theorem 5.1), we have
\vskip-4mm
$$\ba{ll}
\ns\ds \dbE\2n\int_0^T\3n\1n \bar L(t)\big[X_1(t)\1n+\1nX_2(t)\big] dt \1n+\1n\dbE\big[\bar H \big[X_1(T)\1n+\1nX_2(T)\big]\big]\1n=\1n\dbE\2n\int_0^T\3n\1n \lan \f(t), Y(t)\ran \1ndt \1n+\1n\dbE\big[\bar H \f(T)\big].
\ea
$$
\vskip-2mm
\no Let for $j=1,\cds,d$,
\vskip-5.5mm
\bel{decompose}\ba{ll}
\eta(t)\2n:=\2n\begin{pmatrix}\eta^0(t)\\ \eta^1(t)\\ \eta^2(t)\end{pmatrix},
\zeta^j(t)\2n:=\2n\begin{pmatrix}\zeta^{0j}(t)\\ \zeta^{1j}(t)\\ \zeta^{2j}(t)\end{pmatrix},
Y(t)\2n:=\2n\begin{pmatrix}Y^0(t)\\ Y^1(t)\\ Y^2(t)\end{pmatrix},
 Z^j(t,s)\2n:=\2n\begin{pmatrix}Z^{0j}(t,s)\\ Z^{1j}(t,s)\\ Z^{2j}(t,s)\end{pmatrix}.
\ea\ee
\vskip-2mm
\no By \rf{first-order adjoint equation}, we deduce
\vskip-7mm
$$\ba{ll}
\ns\ds \dbE\int_0^T \blan \f(t), Y(t)\bran dt +\dbE\big[\bar H \f(T)\big]\\
\ns\ds =\1n\dbE\2n\int_0^T\3n \int_0^t\2n \blan Y(t),\1n B(t,s)\1n+\1n\bar B(t,s)\bran ds dt \1n+\2n\sum_{j=1}^d\dbE\2n\int_0^T\2n \int_0^t\2n \blan Z^j(t,s), D^j(t,s)\1n+\1n\bar D^j(t,s)\bran ds dt\\
\ns\ds\q +\dbE\bigg[\bar H \int_0^T\big[\bar B(T,s)+B(T,s)\big]ds +\sum_{j=1}^d\int_0^T \zeta^j (s)^\top\big[\bar D^j(T,s)+D^j(T,s)\big]ds\bigg],
\ea
$$
\vskip-4mm
\no which together with \rf{variational-inequality-X} yields that
\bel{variational-inequality-X-mid}\ba{ll}
\ns\ds J(u^\e(\cd))\1n-\1nJ(u^*(\cd))\1n=\1n\dbE\1n\int_0^T\2n \Big[\D l(s)+\frac 1 2 X_1(s)^{\top} L(s) X_1(s)+\Blan\D b(s)\\
\ns\ds \q +\frac 1 2 X_1(s)^\top\pa^2b(s)X_1(s) ,\int_s^T\2nY^0(t)dt+\1n\int_{s+\d}^T\2nY^1(t)dt{\bf1}_{[0,T-\d)}(s)+h_x(T)^\top \\
  \ns\ds \q +h_y(T)^\top{\bf1}_{[0,T-\d)}(s)\Bran+\sum_{j=1}^d\Blan\D\si^j(s)+\frac 1 2 X_1(s)^\top\pa^2\si^j(s)X_1(s)\\
  \ns\ds \q +\D\Xi^j(s)X_1(s),\int_{s+\d}^T\2nZ^{1j}(t,s) dt{\bf1}_{[0,T-\d)}(s)\1n+\2n\int_s^T\3nZ^{0j}(t,s)dt \1n+\1n\zeta^{0j}(s) \1n\\
  \ns\ds \q+\1n\zeta^{1j}(s){\bf1}_{[0,T-\d)}(s)\Bran \Big]ds+\frac 1 2 \dbE X_1(T)^{\top} H X_1(T)+o(\e).
\ea\ee
\vskip-2mm

\ss

Next we would like to write \rf{variational-inequality-X-mid} in a more concise form. To this end, for $j=1,\cds,d$, $0\les\t\les T$, let us denote
\vskip-5.2mm
\bel{pq}\left\{\ba{ll}
\ds \3np(\t)\1n:=\1n\eta^0(\t)\1n+\1n\eta^1(\t){\bf 1}_{[0,T-\d)}(\t) \2n+\1n\dbE_\t\1n\bigg[\1n\int_\t^T\3n\1nY^0(t)dt \1n+\3n\int_{\t+\d}^T\3n\2nY^1(t)dt{\bf1}_{[0,T-\d)}(\t)\bigg],\\
\ns\ds \3nq^j(\t)\1n:=\1n\z^{0j}(\t)\1n+\1n\zeta^{1j}(\t){\bf 1}_{[0,T-\d)}(\t) \1n+\3n\int_\t^T\3n\2nZ^{0j}(t,\t)dt \1n+\3n\int_{\t+\d}^T\3n\2nZ^{1j}(t,\t)dt{\bf1}_{[0,T-\d)}(\t).
\ea\right.\ee

Define $G:[0,T]\times\dbR^n\times\dbR^n\times\dbR^n\times \dbR^{n}\times\dbR^{n\times d}\times\dbR^m\times\dbR^m\ra\dbR$ as
\vskip-1mm
\bel{G}\ba{ll}
\ns\ds G(\t,x,y,z,p,q,u,\m):=l(\t,x,y,z,u,\m)\\
\ns\ds \q+\blan p ,b(\t,x,y,z,u,\m)\bran +\sum_{j=1}^d\blan q^j,\si^j(\t,x,y,z,u,\m)\bran.
\ea\ee

Now we give the main result of this subsection.
\bl\label{lemma variational inequality-1}
Let Assumption {\rm (A1)} hold. Suppose $(x^*(\cdot),u^*(\cdot))$ is an optimal pair, $x^\varepsilon(\cdot)$ is the trajectory corresponding to  $u^\varepsilon(\cdot)$, given by (\ref{perturbed control}), $(\eta(\cd),\z(\cd),Y(\cd),Z(\cd,\cd))$ is a solution to \rf{first-order adjoint equation}. Then, the following variational inequality holds:
\vskip-5mm
\bel{variational inequality-X-1}\ba{ll}
\ns\ds J(u^\e(\cd))\1n-\1nJ(u^*(\cd))\1n=\1n\dbE\2n\int_\t^{\t+\e}\3n\3n\3n\D G(t)dt\1n+\1n\dbE\2n\int_{\t+\d}^{\t+\d+\e}\3n\3n\3n\D\ti{G}(t)dt {\bf1}_{[0,T-\d)}(\t)\1n+\1n\frac 1 2\sE(\e)\1n+\1no(\e),\\
\ea\ee
\vskip-1mm
\no for all $ v(\cd)\in\mathcal{U}_{ad}$ and $\t\in[0,T)$, where
\vskip-1mm
$$\ba{ll}
\ns\ds \D G(t):=G(t,x^*(t),y^*(t),z^*(t),p(t),q(t),v(t),u^*(t-\d))\\
\ns\ds\qq\qq-G(t,x^*(t),y^*(t),z^*(t),p(t),q(t),u^*(t),u^*(t-\d)),
\ea$$
\vskip-1mm
$$\ba{ll}
\ns\ds \D\ti G(t):=G(t,x^*(t),y^*(t),z^*(t),p(t),q(t),u^*(t),v(t-\d))\\
\ns\ds\qq\qq -G(t,x^*(t),y^*(t),z^*(t),p(t),q(t),u^*(t),u^*(t-\d)),
\ea$$
\vskip-1mm
\bel{ee}\ba{ll}
\ns\ds \sE(\e):=\dbE\int_0^T
X_1(t)^\top \pa^2G(t)X_1(t)dt+\dbE \big[X_1(T)^\top HX_1(T)\big].\\
%
%\ns\ds \pa^2G(t):=\begin{bmatrix}
%%
%G_{xx}(t) & G_{xy}(t) & G_{xz}(t) \\
%%
%G_{yx}(t) & G_{yy}(t) & G_{yz}(t) \\
%%
%G_{zx}(t) & G_{zy}(t) & G_{zz}(t)
%%
%\end{bmatrix}.
%
\ea\ee
\el
\proof
\no Notice that
%
%\vskip-7mm
%%
%$$\ba{ll}
%%
%\ns\ds \dbE\int_0^T\Blan\int_s^TZ^{0j}(t,s)dt,\D\si^j(s)\Bran ds\\
%%
%\ns\ds =\dbE\int_\t^{\t+\e}\Blan\int_s^TZ^{0j}(t,s)dt, \si^j(s,x^*(s),y^*(s),z^*(s),v(s),\m^*(s))-\si^j(s,\Th(s))\Bran ds\\
%%
%\ns\ds \q +\dbE\int_{\t+\d}^{\t+\d+\e}\Blan\int_s^T Z^{0j}(t,s)dt,\si^j(s,x^*(s),y^*(s),z^*(s),u^*(s),v(s-\d))\\
%%
%\ns\ds \qq -\si^j(s,\Th(s))\Bran ds {\bf %1}_{[0,T-\d)}(\t),
%%
%\ea$$
%%
%\vskip-2mm
%%
%\no and
%%
$$\ba{ll}
\ns\ds \dbE\int_0^T\2n\int_0^t\2n\Blan Z^{1j}(t,s),{\bf 1}_{(\d,\i)}(t-s)\D\si^j(s)\1n\Bran dsdt \1n=\1n\dbE\int_0^{T-\d}\3n\3n\blan\1n\int_{s+\d}^T\2n Z^{1j}(t,s)dt,\D\si^j(s)\1n\bran ds\\
\ns\ds \1n=\1n\dbE\2n\int_\t^{\t+\e}\3n\2n\blan\1n\int_{s+\d}^T\2n Z^{1j}(t,s)dt,\si^j(s,x^*(s),y^*(s),z^*(s),v(s),\m^*(s)) \1n-\1n\si^j(s,\Th(s))\bran ds \\
\ns\ds \q \times {\bf1}_{[0,T-\d)}(\t)+\dbE\int_{\t+\d}^{\t+\d+\e}\Blan  \si^j\1n(s,x^*(s),y^*(s),z^*(s),u^*(s),v(s-\d))\\
\ns\ds \qq\qq\qq\qq -\si^j(s,\Th(s)),\int_{s+\d}^T Z^{1j}(t,s)dt\Bran ds{\bf1}_{(0,T-\d)}(\t+\d).
\ea$$
\vskip-2mm
\no Then, by applying Lemma \ref{estimate for variational variable}, \rf{variational-inequality-X-mid} and \rf{pq}, we complete the proof.
\endpf

\subsection{Second-order adjoint equations}

To treat the quadratic form in (\ref{variational inequality-X-1}), let us borrow some ideas from \cite{Wang-Yong-2022}.
Now we introduce the following systems of backward equations:
\vskip-6mm
\bel{second-order adjoint equations}\left\{\ba{ll}
\ds (a) \q P_1(r)=H-\sum_{j=1}^d\int_r^T Q_1^j(\th)dW^j(\th),\q 0\les r\les T,\\
\ns\ds (b) \q P_2(r)=A(T,r)^\top P_1(r) +\sum_{j=1}^dC^j(T,r)^\top Q_1^j(r) +\int_r^T\bigg[A(\th,r)^\top P_2(\th)\\
\ns\ds\qq\qq +\sum_{j=1}^dC^j(\th,r)^\top Q_2^j(\th,r)\bigg]d\th\1n-\2n\sum_{j=1}^d\int_r^T\3nQ_2^j(r,\th)dW^j(\th), 0\1n\les\1n r\1n\les\1n T,\\
\ns\ds (c) \q P_3(r)\1n=\1n\pa^2G(r)\1n+\2n\sum_{j=1}^d\1nC^j(T,r)^\top \2nP_1(r)C^j(T,r) \1n\\
\ns\ds\qq\qq +\2n\sum_{j=1}^d\1n\int_r^T\1n\bigg[C^j(T,r)^\top \1nP_2(\th)^\top\1n C^j(\th,r)\1n+C^j\1n(\th,r)^\top\1n P_2(\th)C^j(T,r)\1n\\
\ns\ds \qq\qq+C^j(\th\1n,r)^\top\1n P_3(\th)C^j(\th,r)\1n\bigg]d\th\1n +\3n\int_r^T\1n\3n\int_r^T\3n\1nC^j(\th\1n,r)^\top\2n P_4(\th'\1n,\th)C^j(\th'\1n,r)d\th d\th'\1n\\
\ns\ds \qq\qq-\2n\sum_{j=1}^d\1n\int_r^T\3nQ_3^j(r,\th)dW^j(\th),\q 0\1n\les r\les\1n T,\\
\ns\ds (d)\q  P_4(\th,r)\1n=\1nA(T,r)^\top\1n P_2(\th)^\top \1n +\2n\sum_{j=1}^dC^j(T,r)^\top\1n Q_2^j(\th,r)^\top\1n+A(\th,r)^\top\1n P_3(\th)\\
\ns\ds\qq\qq \1n+\1n\sum_{j=1}^dC^j\1n(\th,r)^\top\1n Q_3^j(\th,r)\1n+\2n\int_r^T\2n\bigg[\sum_{j=1}^dC^j(\th',r)^\top\1n Q_4^j(\th,\th',r)\1n\\
\ns\ds\qq\qq +A(\th',r)^\top\1n P_4(\th,\th')\1n\bigg]d\th'\1n -\2n\sum_{j=1}^d\1n\int_r^T\3nQ_4^j(\th,r,\th')dW^j(\th'), 0\1n\les\1n r\1n\les\1n\th\1n\les\1n T,\\
\ns\ds (e) \q P_4(\th,r)=P_4(r,\th)^\top,\q Q_4(\th,r,\th')=Q_4(r,\th,\th')^\top,\q0\les  \th< r\les T,\\
\ea\right.\ee
\no subject to the following constraints:
\bel{second-order-adjoint-constraints}
\left\{\ba{ll}
\ds P_2(r)=\dbE_\th \big[P_2(r)\big]+\sum_{j=1}^d\int_\th^r Q_2^j(r,\th')dW^j(\th'),\q 0\les r\les T,\\
\ns\ds P_3(r)=\dbE_\th \big[P_3(r)\big]+\sum_{j=1}^d\int_\th^r Q_3^j(r,\th')dW^j(\th'),\q 0\les r\les T,\\
\ns\ds P_4(\th,r)\1n=\1n\dbE_\th' \big[P_4(\th,r)\big]\1n+\1n\sum_{j=1}^d\int_{\th'}^{r\wedge \th}\3nQ_4^j(\th,r,s)dW^j(s), 0\1n\les\1n \th'\1n\les\1n (\th\wedge r) \1n\les\1n T.
\ea\right.
\ee
\no Then, we have the following result for the variational inequality \rf{variational inequality-X-1}.
\bl\label{variational inequality-1}
Let Assumption {\rm (A1)} hold. Suppose $(x^*(\cdot),u^*(\cdot))$ is an optimal pair, $x^\varepsilon(\cdot)$ is the trajectory corresponding to  $u^\varepsilon(\cdot)$, given by (\ref{perturbed control}), $(\eta(\cd),\z(\cd),Y(\cd),Z(\cd,\cd))$ is the solution to \rf{first-order adjoint equation}, $(p(\cd),q(\cd))$ is defined by \rf{pq}. Then, \rf{second-order adjoint equations} admits a unique adapted solution: $(\1nP_1(\cd),\1nQ_1(\cd))\1n\in\1n L^2_{\dbF}(\Omega;C([0,T];\dbS^{3n}\1n))\1n\times\1n \big(L^2_{\dbF}(0,T;\dbS^{3n})\big)^d$, $(P_2(\cd),P_3(\cd),P_4(\cd,\cd))\1n\in \1n L_{\dbF}^2(0,T;\dbR^{(3n)\1n\times \1n(3n)})\1n\times\1n L^2_{\dbF}(0,T;\dbS^{3n})\times L^2\big(0,T;L^2_{\dbF}(0,T;\dbR^{(3n)\times(3n)})\big)$, such that (\ref{second-order-adjoint-constraints}) holds. Furthermore, the variational inequality \rf{variational inequality-X-1} can be deduced as follows:
\vskip-6mm
\bel{variational inequality-X-2}\ba{ll}
\ns\ds J(u^\e(\cd))-J(u^*(\cd))=\dbE\int_\t^{\t+\e}\2n\D G(t)dt+\dbE\int_{\t+\d}^{\t+\d+\e}\2n\D\ti{G}(t)dt{\bf1}_{[0,T-\d)}(\t)\\
\ns\ds \q+\frac 1 2\sum_{j=1}^d\dbE\int_0^T\bigg\{D^j(T,t)^\top P_1(t)D^j(T,t)+\int_t^TD^j(\th,t)^\top P_3(\th)D^j(\th,t)d\th\\
\ns\ds \q+\int_t^T\Big[D^j(T,t)^\top P_2(\th)^\top D^j(\th,t)+D^j(\th,t)^\top P_2(\th)D^j(T,t)\Big]d\th\\
\ns\ds \q+\int_t^T\2n\int_t^TD^j(\th,t)^\top P_4(\th',\th) D^j(\th',t)d\th d\th'\bigg\}dt+o(\e),\q \forall\ \t\in[0,T).
\ea\ee
\vskip-2mm
\el
\proof
\vskip-2mm
\no Note that the BSDE \rf{second-order adjoint equations} (a) admits a unique solution. Then, by Proposition \ref{pro for bsvie} and the similar proof of Theorem 5.1 in \cite{Wang-Yong-2022}, \rf{second-order adjoint equations} has a unique solution.
For simplicity, we just give a sketch of the proof, a detailed proof can be referred to Section 4 in \cite{Wang-Yong-2022}. In the following, without loss of generality, let $d=1$. First we introduce an auxiliary process as follows:
\vskip-5mm
\bel{auxiliary equation}\ba{ll}
\ns\ds \cX_1(t,r)\1n=\1n\int_\t^r\2n\big[A(t,s)X_1(s)\1n+\1nB(t,s)\big]ds \1n+\1n\int_\t^r\1n\big[C(t,s)X_1(s)\1n+\1n D(t,s)\big]dW(s),
\ea\ee
\vskip-1mm
\no for $0\les\t\les r\les t\les T$. Fix $\t=0$, then, $\cX_1(t,t)=X_1(t)$ for all $0\les t\les T$. Applying Lemma \ref{estimate for variational variable}, we have
$\sup\limits_{\t\les t\les T}\dbE\big[\sup\limits_{\t\les r\les t}|\cX_1(t,r)|^p\big]=O(\e^{\frac p 2})$. Let $\Th(\cd,\cd):[\t,T]^2\times \Om\ra\dbR^{(3n)\times (3n)}$ be a process such that for any $t\in[\t,T]$, $\Th(t,\cd)\in L^2_{\dbF}(\t,t;\dbR^{(3n)\times (3n)})$.
Then, by the martingale representation theorem, for any $\t\1n\les\1n s\1n\les\1n t\1n\les\1n T$, there exists a unique $\L(t,s,\cd)\in \big(L^2_{\dbF}(\t,s;\dbR^{(3n)\times (3n)})\big)^d$ satisfying
\vskip-6mm
\bel{Pi-L}\ba{ll}
\ns\ds \Pi(t,s,r)\equiv\dbE_r[\Th(t,s)]=\Th(t,s)-\int_r^s\L(t,s,\th)dW(\th), \ \t\les r\les s\les t\les T.
\ea\ee
\no By It\^o's formula, we obtain for $\t\les r\les s\les t\les T$,
\bel{important relation}\ba{ll}
\ns\ds \dbE\big[\cX_1(t,r)^\top\Th(t,s)\cX_1(s,r)\big]=\dbE\big[\cX_1(t,r)^\top\Pi(t,s,r)\cX_1(s,r)\big]\\
\ns\ds =\dbE\int_\t^r\bigg\{X_1(\th)^\top\[A(t,\th)^\top\Th(t,s) +C(t,\th)^\top\L(t,s,\th)\]\cX_1(s,\th)\\
\ns\ds \q+\cX_1(t,\th)^\top\[\Th(t,s)A(s,\th)+\L(t,s,\th)C(s,\th)\]X_1(\th) +X_1(\th)^\top C(t,\th)^\top\\
\ns\ds \q\times\Th(t,s)C(s,\th)X_1(\th)+D(t,\th)^\top \Th(t,s)D(s,\th)\bigg\}d\th+o(\e).
\ea\ee
\vskip-2mm
\no In the following, we choose different $\Th(\cd,\cd)$, $\Pi(\cd,\cd,\cd)$ and $\L(\cd,\cd,\cd)$ to deal with the quadratic terms about $X_1(\cd)$ in \rf{ee}. First we deal with the term $X_1(T)^\top HX_1(T)$. Take $t=s=T$ and $\Th(T,T)=H$ in \rf{Pi-L}. Then, from \rf{second-order adjoint equations} (a), we have
\vskip-4mm
$$
(\Pi(T,T,r),\L(T,T,r))\equiv(P_1(r),Q_1(r)),\q r\in[0,T].
$$
\vskip-2mm
\no By \rf{important relation}, we get
\vskip-4mm
$$\ba{ll}
\ns\ds \dbE\big[X_1(T)^\top HX_1(T)\big]=\dbE\big[\cX_1(T,T)^\top P_1(T)\cX_1(T,T)\big]\\
\ns\ds =\dbE\int_0^T\bigg\{X_1(r)^\top\[A(T,r)^\top P_1(r)+C(T,r)^\top Q_1(r)\]\cX_1(T,r)\\
\ns\ds \qq+\cX_1(T,r)^\top\[P_1(r)A(T,r)+Q_1(r)C(T,r)\]X_1(r)+X_1(r)^\top C(T,r)^\top\\
\ns\ds \qq  \times P_1(r)C(T,r)X_1(r)+D(T,r)^\top P_1(r)D(T,r)\bigg\}dr+o(\e),
\ea$$
\vskip-1mm
\no which together with \rf{ee} yields that
\vskip-4mm
$$\ba{ll}
\ns\ds \sE(\e)=\dbE\int_0^T\Big\{X_1(r)^\top\[A(T,r)^\top P_1(r)+C(T,r)^\top Q_1(r)\]\cX_1(T,r)\\
\ns\ds \qq\q+\cX_1(T,r)^\top\[P_1(r)A(T,r)+Q_1(r)C(T,r)\]X_1(r)+X_1(r)^\top \[\pa^2G(r)\\
\ns\ds \qq\q +C(T,r)^\top P_1(r)C(T,r)\]X_1(r)+D(T,r)^\top P_1(r)D(T,r)\Big\}dr +o(\e).
\ea$$
\vskip-1mm

Next we deal with the term $X_1(r)^\top[\cds]\cX_1(T,r)$ and $\cX_1(T,r)^\top[\cds]X_1(r)$. Take $t=T$ in \rf{Pi-L}, let
\vskip-4mm
$$
(\Th(T,r),\L(T,\th,r))\equiv(P_2(r)^\top,Q_2(\th,r)^\top),\q 0\les r\les \th\les T.
$$
\vskip-1mm
\no Then, by \rf{important relation} and \rf{second-order adjoint equations} we obtain
\vskip-6mm
\bel{ee-2}\ba{ll}
\ns\ds \sE(\e)=\dbE\int_0^T\bigg\{X_1(r)^\top\bigg[\pa^2G(r)+C(T,r)^\top P_1(r)C(T,r)\\
\ns\ds \qq\q+\int_r^T\(C(T,r)^\top P_2(\th)^\top C(\th,r)+C(\th,r)^\top P_2(\th)C(T,r)\)d\th\bigg]X_1(r)\\
\ns\ds \q\qq+\int_r^T\bigg[X_1(r)^\top\(P_2(\th)A(T,r) +Q_2(\th,r)C(T,r)\)^\top\cX_1(\th,r)\\
\ns\ds \qq\q+\cX_1(\th,r)^\top \(P_2(\th)A(T,r) +Q_2(\th,r)C(T,r)\) X_1(r)\bigg] d\th\bigg\}dr\\
\ns\ds \qq\q+\dbE\int_0^T\bigg\{\int_r^T\[D(T,r)^\top P_2(\th)^\top D(\th,r)+D(\th,r)^\top P_2(\th)D(T,r)\]d\th\\
\ns\ds \qq\q+D(T,r)^\top P_1(r)D(T,r)\bigg\}dr+o(\e).
\ea\ee
\vskip-2mm

Finally we eliminate the terms $X_1(r)^\top[\cds]X_1(r)$, $\cX_1(\th,r)^\top[\cds]X_1(r)$ and their transpose. Take $t=s$ in \rf{Pi-L} and let
\vskip-4mm
$$\ba{ll}
\Th(\th,\th)\equiv P_3(\th),\q \L(\th,\th,r)\equiv Q_3(\th,r),\q 0\les r\les \th\les T.
\ea$$
\vskip-1mm
\no Then, from \rf{important relation} we derive
\vskip-6mm
\bel{ee-3}\ba{ll}
\ns\ds \dbE\int_0^TX_1(r)^\top \Th(r,r)X_1(r)dr\\
\ns\ds =\1no(\e)\1n+\1n\dbE\int_0^T\3n\int_r^T\3n\Big\{X_1(r)^\top\1n \[A(\th,r)^\top\1n \Th(\th,\th)\1n+\1nC(\th,r)^\top\L(\th,\th,r)\]\1n\cX_1(\th,r)\\
\ns\ds \qq+\cX_1(\th,r)^\top\1n\[A(\th,r)^\top \Th(\th,\th)^\top\1n+\1nC(\th,r)^\top\1n\L(\th,\th,r)^\top\1n\]^\top\2n X_1(r)\1n\\
\ns\ds \qq+\1nX_1(r)^\top\1n C(\th,r)^\top\Th(\th,\th)C(\th,r)X_1(r) \1n+\1nD(\th,r)^\top\1n\Th(\th,\th)D(\th,r)\1n\Big\}d\th dr\1n.
\ea\ee
\vskip-3mm
\no Let
\vskip-6mm
$$\ba{ll}
\Th(\th,\th')=P_4(\th,\th')^\top,\q \L(\th,r,\th')=Q_4(\th,r,\th')^\top,\q 0\les\th'\les r\les \th\les T.
\ea$$
\vskip-1mm
\no Then, by \rf{important relation} we get
\vskip-4mm
$$\ba{ll}
\ns\ds \dbE\2n\int_0^T\2n\int_r^T\2n\cX_1(\th,r)^\top\1n\Th(\th,r)X_1(r)d\th dr\1n=\1n\dbE\1n\int_0^T\2n\bigg\{\2n\int_r^T\2n\int_\th^T \2nX_1(r)^\top\1n\[A(\th',r)^\top\1n\Th(\th',\th)\\
\ns\ds \q+C(\th',r)^\top\L(\th',\th,r)\]\cX_1(\th,r)d\th'd\th +\int_r^T\int_r^\th\cX_1(\th,r)^\top\[A(\th',r)^\top\Th(\th,\th')^\top  \\
 \ns\ds\q +C(\th',r)^\top\L(\th,\th',r)^\top\]^\top X_1(r)d\th'd\th+\int_r^T\int_\th^T\[X_1(r)^\top C(\th',r)^\top\\
 \ns\ds \q \times\Th(\th',\th)C(\th,r)X_1(r)+D(\th',r)^\top\Th(\th',\th)D(\th,r)\]d\th'd\th \bigg\}dr+o(\e),
\ea$$
\vskip-2mm
\no which and \rf{variational inequality-X-1}, \rf{ee-2}, \rf{ee-3} imply that \rf{variational inequality-X-2} holds.
\endpf

\br
It is worth mentioning that the first-order adjoint equation \rf{first-order adjoint equation}, consisting of a BSDE and a backward SVIE, is dual with the first-order and second-order variational equations \rf{variational equation-1}-\rf{variational equation-2}, and the second-order adjoint equation \rf{second-order adjoint equations}, consisting of a BSDE and three coupled backward SVIEs, is dual with $(x_1(t)^\top\2n,y_1(t)^\top,\2nz_1(t)^\top\1n)[\cds]
(x_1(t)^\top\2n,y_1(t)^\top\2n,z_1(t)^\top)^\top$, even though the pointwise state delay appears in the state equation and the terminal cost.
\er

\br\label{Remark4.4}
To deal with the cross term $x_1(t)^\top[\cds]y_1(t)$ and its transpose, \cite{Meng-Shi-2021} introduced a new BSDE but required its solution to be zero. In this paper, we get rid of this strict condition. First the delayed variational equations \rf{variational equation-1}-\rf{variational equation-2} are transformed into the Volterra integral equations without delay \rf{variational-1-X}-\rf{variational-2-X}, so that the delayed finite dimensional control problem is converted into another finite dimensional control problem without delay. Then from the above proof, $X_1(r)^\top[\cds]X_1(r)$ contains the cross terms $x_1(t)^\top[\cds]y_1(t)$ and $y_1(t)^\top[\cds]x_1(t)$, so the auxiliary equation \rf{auxiliary equation} is constructed and the set of backward SVIEs \rf{second-order adjoint equations} is introduced to deal with the ``cross terms", without any additional conditions.
\er

\section{General maximum principle}

In this section, we obtain a general maximum principle for Problem (P), and further express first-order and second-order adjoint equations in more compact forms.

\subsection{General maximum principle}

First let us do some interesting analysis of the second-order adjoint equation \rf{second-order adjoint equations}. In the following, we suppose $\t\in[0,T)$ and define
\vskip-4mm
$$\ba{ll}
\ns\ds P_k(\cd):=\left\{\begin{array}{ccccc}
P_k^{(11)}(\cd) &  P_k^{(12)}(\cd) & P_k^{(13)}(\cd)  \\
P_k^{(21)}(\cd) & P_k^{(22)}(\cd) & P_k^{(23)}(\cd)\\
P_k^{(31)}(\cd) & P_k^{(32)}(\cd) & P_k^{(33)}(\cd)
\end{array}
\right\},\ \  k=1,2,3,4.
\ea$$
\vskip-3mm

\ms

\textbf{Case I: The term of $(P_1,Q_1)$.}

By the definition of $H$, we see that
\vskip-4mm
\bel{P1}\ba{ll}
\ns\ds P_1^{(i\ell)}(r)=h_{\k^i\k^\ell}(T) -\sum_{j=1}^d\int_r^TQ_{1j}^{(i\ell)}(\th)dW^j(\th),\q \t\les r\les T,
\ea\ee
\vskip-2mm
\no where $i,\ell=1,2,3$, and $\k^1:=x$, $\k^2:=y$, $\k^3:=z$. In addition,
\vskip-5mm
\bel{DDP1}\ba{ll}
\ns\ds D^j(T,\t)^\top P_1(\t)D^j(T,\t)=\D \si^j (\t)^\top P_1^{(11)}(\t)\D \si^j(\t)\\
\ns\ds \qq +\D\si^j(\t)^\top\big[P_1^{(12)}(\t)+P_1^{(21)}(\t) +P_1^{(22)}(\t)\big]\D\si^j(\t){\bf1}_{(\d,\i)}(T-\t).
\ea\ee
\vskip-3mm

\ms

\textbf{Case II: The term of $(P_2,Q_2)$.}

Let us look at $(P_2,Q_2)$ in \rf{second-order adjoint equations},
\vskip-8mm
$$\ba{ll}
\ns\ds P_2^{(i\ell)}(r)=\psi_2^{(i\ell)}(r)+\int_r^T g_2^{(i\ell)}(\th,r)d\th-\sum_{j=1}^d\int_r^TQ_{2j}^{(i\ell)}(r,\th)dW^j(\th),\q \t\les r\les T,
\ea$$
\vskip-4mm
\no where $i,\ell=1,2,3$. Set
\vskip-6mm
$$\ba{ll}
\ns\ds \Big\{g_2^{(i\ell)}(\th,r)\Big\}_{i,\ell=1}^3:=A(\th,r)^\top P_2(\th)+\sum_{j=1}^dC^j(\th,r)^\top Q_{2}^j(\th,r),
\ea$$
\vskip-4mm
$$\ba{ll}
\ns\ds \Big\{\psi^{(i\ell)}_2(r)\Big\}_{i,\ell=1}^3:=A(T,r)^\top P_1(r)+\sum_{j=1}^dC^j(T,r)^\top Q_1^j(r).
\ea$$
\vskip-2mm

\no For $j=1,\cds,d$,\ $\ell=1,2,3$ and $\k^1:=x$, $\k^2:=y$, $\k^3:=z$, define for $\t\les r\les T$,
\vskip-7mm
$$\ba{ll}
\ns\ds \cG_2^{(\ell)}(r):=h_{x\k^\ell}(T)+  \int_r^T P_2^{(1\ell)}(\th)d\th+\bigg[h_{y\k^\ell}(T)+ \int_{r+\d}^T P_2^{(2\ell)}(\th)d\th\bigg]{\bf1}_{[0,T-\d)}(r),\\
\ns\ds \cQ_{2j}^{(\ell)}(r)\1n:=\1nQ_{1j}^{(1\ell)}(r)\1n+\1n\int_r^T \3nQ_{2j}^{(1\ell)}(\th,r)d\th\1n +\1n\bigg[Q_{1j}^{(2\ell)}(r)\1n+\1n\int_{r+\d}^T \3n Q_{2j}^{(2\ell)}(\th,r)d\th\bigg]{\bf1}_{[0,T-\d)}(r),\\
\ns\ds \cK_2^{(\ell)}(r):=P_1^{(3\ell)}(r) +\int_r^TP_2^{(3\ell)}(\th)d\th.
\ea$$
\vskip-2mm
\no Then, we deduce that for $\t\les r\les T$,
\vskip-5mm
\bel{P2-1}\ba{ll}
\ns\ds P_2^{(1\ell)}(r)=\dbE_r\bigg[b_x(r)^\top\cG_2^{(\ell)}(r) +\sum_{j=1}^d\si_x^j(r)^\top\cQ_{2j}^{(\ell)}(r) +\cK_2^{(\ell)}(r)\bigg],\\
\ns\ds P_2^{(2\ell)}(r)\1n=\1n\dbE_r\bigg[b_y(r)^\top\1n\cG_2^{(\ell)}(r) \1n+\2n\sum_{j=1}^d\1n\si_y^j(r)^\top\1n\cQ_{2j}^{(\ell)}(r) \1n-\1ne^{-\l\d}\cK_2^{(\ell)}(r)\bigg],\\
\ns\ds P_2^{(3\ell)}(r)=\dbE_r\bigg[b_z(r)^\top\cG_2^{(\ell)}(r) +\sum_{j=1}^d\si_z^j(r)^\top\cQ_{2j}^{(\ell)}(r)-\l\cK_2^{(\ell)}(r)\bigg].
\ea\ee
\no For the $P_2$ part in (\ref{variational inequality-X-2}), we have
\vskip-4mm
\bel{DDP2}\ba{ll}
\ns\ds \dbE_\t\int_\t^T\Big[D^j(T,\t)^\top P_2(\th)^\top D^j(\th,\t)+D^j(\th,\t)^\top P_2(\th)D^j(T,\t)\Big]d\th\\
\ns\ds =\D\si^j(\t)^\top\dbE_\t\bigg[\int_\t^T\(P_2^{(11)}(\th)^\top +P_2^{(11)}(\th) +\big[P_2^{(12)}(\th)^\top\\
\ns\ds \q+P_2^{(12)}(\th)\big]{\bf1}_{[0,T-\d)}(\t)\)d\th +\int_{\t+\d}^T \[P_2^{(21)}(\th)^\top+P_2^{(21)}(\th)\\
\ns\ds \q +P_2^{(22)}(\th)^\top+P_2^{(22)}(\th)\]d\th{\bf1}_{[0,T-\d)}(\t)\bigg]\D\si^j(\t).
\ea\ee

\ms

\textbf{Case III: The term of $(P_4,Q_4)$.}

Let us look at $(P_4,Q_4)$ in \rf{second-order adjoint equations},
\vskip-4mm
$$\ba{ll}
\ns\ds  P_4^{(i\ell)}\1n(\th,r)\1n=\1n\psi^{(i\ell)}_4\1n(\th,r) \1n+\1n\int_r^T\3n g^{(i\ell)}_4(\th,\th',r)d\th'\1n-\2n\sum_{j=1}^d \int_r^T\3nQ_{4j}^{(i\ell)}(\th,r,\th')dW^j(\th'),
\ea$$
\vskip-1mm
\no where $\t\les r\les\th \les T$, $i,\ell=1,2,3$. Define
\vskip-4mm
\bel{psi4-ij}\ba{ll}
\ns\ds \Big\{\psi^{(i\ell)}_4(\th,r)\Big\}_{i,\ell=1}^3:=\1nA(T,r)^\top P_2(\th)^\top +\1n\sum_{j=1}^dC^j(T,r)^\top Q_2^j(\th,r)^\top \\
\ns\ds \qq\qq\qq\qq\q +A(\th,r)^\top P_3(\th)+\sum_{j=1}^dC^j(\th,r)^\top Q_3^j(\th,r),
\ea\ee
\vskip-2mm
\no and
\vskip-5mm
\bel{g4-ij}\ba{ll}
\ns\ds \Big\{g_4^{(i\ell)}(\th,\th',r)\Big\}_{i,j=1}^3 :=A(\th',r)^\top P_4(\th,\th')+\sum_{j=1}^dC^j(\th',r)^\top Q_4^j(\th,\th',r).
\ea\ee
\vskip-1mm
\no For $\ell=1,2,3$, $j=1,\cds,d$ and $\th\ges r$, define
\vskip-3mm
\bel{cG3-CQ3}\ba{ll}
\ns\ds \cG_4^{(\ell)}(\th,r):=P_2^{(\ell1)}(\th)^\top +P_3^{(1\ell)}(\th) +{\bf1}_{(\d,\infty)}(\th-r) P_3^{(2\ell)}(\th)+\1n{\bf1}_{(\d,\i)}(T\1n-\1nr) \\
\ns\ds \qq\qq \times P_2^{(\ell2)}(\th)^\top\1n+\1n\int_r^T \3nP_4^{(1\ell)}(\th,\th')d\th'\1n +\1n{\bf1}_{(\d,\infty)}(T\1n-\1nr)\1n\int_{r+\d}^T\3nP_4^{(2\ell)}(\th,\th')d\th',\\
\ns\ds \cQ_{4j}^{(\ell)}(\th,r):=Q_{2j}^{(\ell1)}(\th,r)^\top +Q_{3j}^{(1\ell)}(\th,r) +{\bf1}_{(\d,\infty)}(\th-r) Q_{3j}^{(2\ell)}(\th,r) \\
\ns\ds\qq\qq+{\bf1}_{(\d,\i)}(T-r)Q_{2j}^{(\ell2)}(\th,r)^\top +\int_r^T Q_{4j}^{(1\ell)}(\th,\th',r)d\th'\\
\ns\ds \qq\qq+\int_{r+\d}^TQ_{4j}^{(2\ell)}(\th,\th',r)d\th'{\bf1}_{(\d,+\infty)}(T-r),\\
\ns\ds \cK_4^{(\ell)}(\th,r):=P_2^{(\ell3)}(\th)^\top +P_3^{(3\ell)}(\th)+\int_r^TP_4^{3\ell}(\th,\th')d\th'.
\ea\ee
\vskip-1mm
\no Then, for $\th\ges r$, we have
\vskip-5mm
\bel{P3-1}\ba{ll}
\ns\ds P_4^{(1\ell)}(\th,r)=\dbE_r \Big[b_x(r)^\top \cG_4^{(\ell)}(\th,r)+\sum_{j=1}^d\si_x^j(r)^\top \cQ_{4j}^{(\ell)}(\th,r)+\cK_4^{(\ell)}(\th,r)\Big],\\
\ns\ds P_4^{(2\ell)}(\th,r)=\dbE_r \Big[b_y(r)^\top \cG_4^{(\ell)}(\th,r)+\sum_{j=1}^d\si_y^j(r)^\top \cQ_{4j}^{(\ell)}(\th,r)-e^{-\l\d}\cK_4^{(\ell)}(\th,r)\Big],\\
\ns\ds P_4^{(3\ell)}(\th,r)=\dbE_r \Big[b_z(r)^\top \cG_4^{(\ell)}(\th,r)+\sum_{j=1}^d\si_z^j(r)^\top \cQ_{4j}^{(\ell)}(\th,r)-\l\cK_4^{(\ell)}(\th,r)\Big].
\ea\ee
\vskip-3mm
\no For $\th< r$, set
\vskip-4mm
$$\ba{ll}
\ns\ds P_4^{(i\ell)}(\th,r):=P_4^{(\ell i)}(r,\th)^\top,\q Q_4^{(i\ell)}(\th,\th',r):=Q_4^{(\ell i)}(\th',\th,r)^{\top},\q i,\ell=1,2,3.
\ea$$
\vskip-1mm
\no Next we look at the $P_4$ part in (\ref{variational inequality-X-2}). Denote
\vskip-6mm
$$\ba{ll}
\ns\ds \cP_4(\t):=\int_\t^T\int_\t^T P_4^{(11)}(\th',\th) d\th d\th' +\bigg(\int_{\t+\d}^T \int_\t^T P_4^{(12)}(\th',\th) d\th d\th'\\
\ns\ds \qq\qq +\int_\t^T \int_{\t+\d}^T P_4^{(21)}(\th',\th) d\th d\th' +\int_{\t+\d}^T \int_{\t+\d}^T P_4^{(22)}(\th',\th) d\th d\th'\bigg){\bf1}_{[0,T-\d)}(\t).
\ea$$
\vskip-3mm
\no Then, we have
\vskip-4mm
\bel{DDP3}\ba{ll}
\ns\ds \dbE_\t\1n\bigg[\1n\int_\t^T\1n\int_\t^T\3nD^j(\th,\t)^\top\1n P_4(\th',\th) D^j(\th',\t)d\th d\th'\1n\bigg] \1n=\1n\D \si^j(\t)^\top\1n\dbE_\t\big[\cP_4(\t)\big]\D\si^j(\t).
\ea\ee
\vskip-1mm
%
%\vskip-8mm
%%
%\bel{P3-2}\ba{ll}
%%
%\ns\ds P_4^{(\ell1)}(\th,r)=\dbE_\th \bigg[\cG_4^{(\ell)}(r,\th)^\top b_x(\th) +\sum_{j=1}^d\cQ_{4j}^{(\ell)}(r,\th)^\top\si_x^j(\th) +\cK_4^{(\ell)}(r,\th)^\top\bigg],\\
%%
%\ns\ds P_4^{(\ell2)}(\th,r)=\dbE_\th \bigg[\cG_4^{(\ell)}(r,\th)^\top b_y(\th) +\sum_{j=1}^d\cQ_{4j}^{(\ell)}(r,\th)^\top\si_y^j(\th) -e^{-\l\d}\cK_4^{(\ell)}(r,\th)^\top\bigg],\\
%%
%\ns\ds P_4^{(\ell3)}(\th,r)=\dbE_\th \bigg[\cG_4^{(\ell)}(r,\th)^\top b_z(\th) +\sum_{j=1}^d\cQ_{4j}^{(\ell)}(r,\th)^\top\si_z^j(\th) -\l\cK_4^{(\ell)}(r,\th)^\top\bigg].
%%
%\ea\ee
%%
%\vskip-4mm

\ms

\textbf{Case IV: The term of $(P_3,Q_3)$.}

Now, let us look at $(P_3,Q_3)$ in (\ref{second-order adjoint equations}),
\vskip-4mm
$$\ba{ll}
\ns\ds P_3^{(i\ell)}(r)=\psi_3^{(i\ell)}(r)+\int_r^T g_3^{(i\ell)}(\th,r)d\th-\sum_{j=1}^d\int_r^TQ_{3j}^{(i\ell)}(r,\th)dW^j(\th),\\
\ns\ds \qq\qq\qq\qq\qq \t\les r\les T,\q i,\ell=1,2,3.
\ea$$
\vskip-2.5mm
\no Define
\vskip-8.5mm
$$\ba{ll}
\ns\ds \Big\{\1n\psi_3^{(i\ell)}(r)\1n\Big\}_{i,\ell=1}^3\3n:=\1n\pa^2G(r) \2n+\2n\sum_{j=1}^d \1nC^j(T,r)^\top\1n P_1(r)C^j(T,r)\1n+\2n\sum_{j=1}^d\1n\int_r^T\3n\3n\big[C^j(T,r)^\top\1n P_2(\th)^\top\1n C^j(\th,r)\\
\ns\ds\qq\qq+C^j(\th,r)^\top P_2(\th)C^j(T,r)\big]d\th \1n+\2n\sum_{j=1}^d \int_r^T\3n\int_r^T\3n C^j(\th,r)^\top P_4(\th',\th)C^j(\th',r) d\th d\th',
\ea$$
\vskip-3mm
\no and
\vskip-6mm
$$\ba{ll}
\ns\ds \Big\{g_3^{(i\ell)}(\th,r)\Big\}_{i,\ell=1}^3:=\sum_{j=1}^d C^j(\th,r)^\top P_3(\th)C^j(\th,r).
\ea$$
\vskip-2mm
\no Then, for $\ell=1,2,3$, and $\k^1:=x$, $\k^2:=y$, $\k^3:=z$, we have
\vskip-4mm
\bel{P4}\ba{ll}
\ns\ds P_3^{(1\ell)}(r)=G_{x\k^\ell}(r)+\sum_{j=1}^d\si_{x}^j(r)^\top \dbE_r \big[\cP(r)\big]\si_{\k^\ell}^j(r),\q \t\les r\les T,\\
\ns\ds P_3^{(2\ell)}(r)=G_{y\k^\ell}(r)+\sum_{j=1}^d\si_{y}^j(r)^\top\dbE_r \big[\cP(r)\big]\si_{\k^\ell}^j(r),\q \t\les r\les T,\\
\ns\ds P_3^{(3\ell)}(r)=G_{z\k^\ell}(r)+\sum_{j=1}^d\si_z^j(r)^\top\dbE_r \big[\cP(r)\big]\si_{\k^\ell}^j(r),\q \t\les r\les T,
\ea\ee
\no where
$$\ba{ll}
\ns\ds \cP(r):=h_{xx}(T)+\big[h_{yx}(T)+h_{xy}(T)+h_{yy}(T)\big]{\bf1}_{[0,T-\d)}(r) \\
\ns\ds \q +\int_r^T\3n\[P_2^{(11)}(\th)^\top\1n+\1nP_2^{(11)}(\th)\]d\th\1n +\1n\int_r^T\[P_2^{(12)}(\th)^\top+P_2^{12}(\th)\]d\th{\bf1}_{[0,T-\d)}(r)\\
\ns\ds \q +\int_{r+\d}^T\[P_2^{(21)}(\th)^\top+ P_2^{(21)}(\th)+P_2^{(22)}(\th)^\top+P_2^{(22)}(\th)\]d\th{\bf1}_{[0,T-\d)}(r)\\
\ns\ds\q +\int_r^T\int_r^TP_4^{(11)}(\th',\th) d\th d\th' +\bigg\{\int_{r+\d}^T \int_r^T P_4^{(12)}(\th',\th) d\th d\th'\\
\ea$$
\bel{cP}\ba{ll}
\ns\ds\q +\int_r^T \int_{r+\d}^T P_4^{(21)}(\th',\th) d\th d\th' +\int_{r+\d}^T \int_{r+\d}^T P_4^{(22)}(\th',\th) d\th d\th'\bigg\}{\bf1}_{[0,T-\d)}(r)\\
\ns\ds\q  +\int_r^T P_3^{(11)}(\th)d\th+\int_{r+\d}^T\[P_3^{(21)}(\th)+P_3^{(12)}(\th)+P_3^{(22)}(\th)\]d\th{\bf1}_{[0,T-\d)}(r)\\
\ns\ds  =\cG_2^{(1)}(r)+\int_r^T \cG^{(1)}_4(\th',r)d\th'+\bigg[\cG_2^{(2)}(r)+\int_{r+\d}^T \cG^{(2)}_4(\th',r)d\th'\bigg]{\bf1}_{[0,T-\d)}(r).
\ea\ee
\vskip-1mm
\no Next, for the $P_3$ part in (\ref{variational inequality-X-2}), we have
\vskip-3mm
\bel{DDP4}\ba{ll}
\ns\ds \dbE_\t\bigg[\int_\t^T D^j(\th,\t)^\top P_3(\th)D^j(\th,\t)d\th\bigg]=\D \si^j(\t)^\top \dbE_\t\bigg\{\int_\t^T P_3^{(11)}(\th) d\th\\
\ns\ds \qq +\int_{\t+\d}^T \[P_3^{(12)}(\th)+P_3^{(21)}(\th)+P_3^{(22)}(\th)\]d\th{\bf1}_{(\d,\i)}(T-\t)\bigg\}\D \si^j(\t).
\ea\ee
\vskip-1mm

\ss

Based on the above preparation, now we are in a position to state the general maximum principle for Problem (P). Recall \rf{G} and define the Hamiltonian function $\cH:[0,T]\times\dbR^n\times\dbR^n\times\dbR^n \times\dbR^n\times\dbR^{n\times d}\times\dbS^n\times\dbR^m\times\dbR^m\ra\dbR$ as follows:
\vskip-8mm
$$\ba{ll}
\ns\ds \cH(\t,x,y,z,p,q,\cP,u,\m) :=G(\t,x,y,z,p,q,u,\m)+\sum_{j=1}^dTr\[\big(\si^j(\t,x,y,z,u,\m) \\
\ns\ds \qq\qq\qq\qq\qq\q  -\si^j(\t,\Th(\t))\big)^\top\cP\big(\si^j(\t,x,y,z,u,\m) -\si^j(\t,\Th(\t))\big)\].
\ea$$

\bt\label{MP}
Let Assumption {\rm (A1)} hold. Suppose $(x^*(\cdot),u^*(\cdot))$ is an optimal pair, $(\eta(\cd),\z(\cd),Y(\cd),\\Z(\cd,\cd))$ is the solution to \rf{first-order adjoint equation}, $(p(\cd),q(\cd))$ and $\cP(\cd)$ are defined by \rf{pq} and \rf{cP}, $(P_1(\cd),P_2(\cd),P_3(\cd),\\P_4(\cd,\cd))$ is the solution to \rf{second-order adjoint equations}-\rf{second-order-adjoint-constraints}. Then, the following maximum condition holds:
\vskip-3mm
\bel{maximum condition}\ba{ll}
\ns\ds \D\cH(\t)+\dbE_\t\big[\D\ti \cH(\t+\d){\bf1}_{[0,T-\d)}(\t)\big]\ges0, \q \forall\ v\in U,\q \ae,\ \as
\ea\ee
\vskip-1mm
\no where
\vskip-3mm
$$\ba{ll}
\ns\ds \D\cH(\t):=\cH(\t,x^*(\t),y^*(\t),z^*(\t),p(\t),q(\t), \cP(\t),v,\m^*(\t))\\
\qq\qq-\cH(\t,x^*(\t),y^*(\t),z^*(\t),p(\t),q(\t), \cP(\t),u^*(\t),\m^*(\t)),
\ea$$
\vskip-1mm
$$\ba{ll}
\ns\ds \D\ti\cH(\t):=\cH(\t,x^*(\t),y^*(\t),z^*(\t),p(\t),q(\t), \cP(\t),u^*(\t),v)\\
\ns\ds \qq\qq-\cH(\t,x^*(\t),y^*(\t),z^*(\t),p(\t),q(\t), \cP(\t),u^*(\t),\m^*(\t)).
\ea$$

\et
\proof
By Lemma \ref{variational inequality-1}, \rf{DDP1}, \rf{DDP2}, \rf{DDP3}, \rf{DDP4} and \rf{cP}, we obtain
$$\ba{ll}
\ns\ds J(u^\e(\cd))-J(u^*(\cd))=\dbE\int_\t^{\t+\e}\D G(t)dt+\dbE\int_{\t+\d}^{\t+\d+\e}\D\ti{G}(t)dt{\bf1}_{[0,T-\d)}(\t)\\
\ns\ds \qq\qq\qq\qq\q +\frac 1 2 \sum_{j=1}^d\dbE\int_0^TTr[\D\si^j(t)^\top\cP(t)\D\si^j(t)]dt+o(\e).
\ea$$
\vskip-1mm
\no Thus, similar to the proof of Theorem 4.1 in \cite{Meng-Shi-2021}, we complete the proof.
\endpf
\br
The general maximum principle \rf{maximum condition} consists of two parts: $\dbE_{\t}[\D\ti {\cH}(\t+\d)]$ characterizes the maximum condition with delay, while $\D\cH(\t)$ characterizes the one without delay, in similar form to (3.20) in Chapter 3 of \cite{Yong-Zhou-1999}.
\er

\br
Compared with \cite{Meng-Shi-2021}, (i) when the distributed delay appears in the control system, the general maximum principle of optimal control for stochastic differential delay systems can be obtained; (ii) the maximum condition \rf{maximum condition} is similar to (5.6) in \cite{Meng-Shi-2021}, but all the additional requirements in \cite{Meng-Shi-2021} are removed; (iii) a new set of backward SVIEs \rf{second-order adjoint equations} is introduced to deal with the ``cross term", instead of the special BSDE (5.3) in \cite{Meng-Shi-2021}.
\er

%\br
%%
%\vskip-1.5mm
%%
%Noting $(\eta^2(\cd),\zeta^2(\cd),Y^2(\cd),Z^2(\cd,\cd))$ and $P_k^{(3\ell)}$, $P_k^{(\ell 3)}$ do not explicitly appear in the maximum condition \rf{maximum condition}, $k=1,2,3,4$, $\ell=1,2,3$. In fact, them are implicitly contained in \rf{maximum condition} from the coupled structure of \rf{first-order adjoint equation} and \rf{P2-1}, \rf{P3-1}, \rf{P3-2}, \rf{P4}.
%%
%\er

\subsection{Extensions of adjoint equations}

\vskip-1mm

In this subsection, we further explore the first-order and second-order adjoint equations \rf{first-order adjoint equation} and \rf{second-order adjoint equations}. Interestingly, under some cases, \rf{first-order adjoint equation} and \rf{second-order adjoint equations} have more compact forms, similar to the existing literatures \cite{Meng-Shi-2021,Yong-Zhou-1999,Zhang-2022}.

\subsubsection{Extensions of first-order adjoint equations}

\vskip-1mm

We rewrite \rf{pq}, and define $(\ti{p}(\cd),\ti{q}(\cd))$ as follows: for $j=1,\cds,d$,
\vskip-5mm
\bel{pq-}\left\{\ba{ll}
\ds\3n p(\t)\1n:=\1n\eta^0(\t)\1n+\1n\eta^1(\t){\bf 1}_{[0,T-\d)}(\t)\1n+\1n \dbE_\t\bigg[\1n\int_\t^T\3nY^0(t)dt \1n+\1n\int_{\t+\d}^T\3n\3nY^1(t)dt{\bf1}_{[0,T-\d)}(\t)\bigg],\\
\ns\ds \3n q^j(\t)\1n:=\1n\z^{0j}\1n(\t)\1n+\1n\zeta^{1j}(\t){\bf 1}_{[0,T\1n-\d)}(\t)\1n +\3n\int_\t^T\3n\1nZ^{0j}(t,\1n\t)dt \1n+\3n\int_{\t+\d}^T\3n\2nZ^{1j}(t,\1n\t)dt{\bf1}_{[0,T\1n-\d)}(\t),\\
\ns\ds \3n \ti{p}(\t)\1n:=\1n\dbE_\t\bigg[\int_\t^T\3nY^2(t)dt\bigg]\1n+\1n\eta^2(\t),\q \ti{q}^j(\t)\1n:=\1n\int_\t^T\3nZ^{2j}(t,\t)dt\1n+\1n\zeta^{2j}(\t).
\ea\right.\ee
\vskip-2mm

Now we can link the first-order adjoint equation \rf{first-order adjoint equation} with a set of anticipated BSDEs.

\bt\label{theorem relation-1}
Let Assumption {\rm (A1)} hold. Suppose $(x^*(\cdot),u^*(\cdot))$ is an optimal pair, $(\eta(\cd),\z(\cd),Y(\cd),\\Z(\cd,\cd))$ is the solution to \rf{first-order adjoint equation}. Then, $(p(\cd),q(\cd),\ti p(\cd),\ti q(\cd))$ defined by \rf{pq-} satisfies the following set of anticipated BSDEs:
\vskip-6mm
\bel{pq- equation}\left\{\ba{ll}
\ds \3n p(t)\1n=\1nh_x(T)^\top\2n+\1n\int_t^T\3n\bigg\{b_x(s)^\top p(s)\1n+\2n\sum_{j=1}^{d}\1n\si_x^j(s)^\top\1n q^j(s)\1n+l_x(s)^\top\1n+\1n \ti p(s)\1n\bigg\}ds\\
 \ns\ds \qq\q \1n-\2n\sum_{j=1}^d\1n\int_t^T\3nq^j(s)dW^j(s),\q t\in[T-\d,T],\\
\ns\ds \3n p(t)\1n=\1np(T-\d)\2n+\1n\dbE_{T-\d}\big[h_y(T)^\top \big] \1n+\1n\int_t^{T-\d}\2n\bigg\{b_x(s)^\top p(s)\1n+\2n\sum_{j=1}^{d}\1n\si_x^j(s)^\top\1n q^j(s)\1n\\
\ns\ds \qq\q +l_x(s)^\top\1n+\1n \ti p(s)\1n+\1n\dbE_s\big[b_y(s\1n+\1n\d)^\top\1n p(s\1n+\1n\d)\1n+\2n\sum_{j=1}^d\si_y^j(s+\d)^\top \1nq^j(s\1n+\1n\d)\1n \\
\ns\ds \qq\q+l_y(s\1n+\1n\d)\1n^\top\3n-\1ne^{-\l\d}\ti p(s\1n+\1n\d)\big] \1n\bigg\}ds \1n-\3n\sum_{j=1}^d\1n\int_t^{T-\d}\3n\1nq^j(s)dW^j(s), t\1n\in\1n[0,\1nT\1n-\1n\d),\\
\ns\ds \3n \ti p(t)\1n=\1nh_z(T)^\top\1n+\3n\int_t^T\3n\Big\{b_z(s)^\top\1n p(s)\1n+\2n\sum_{j=1}^d\si_z^j(s)^\top\1n q^j(s)\1n+\1nl_z(s)^\top\1n-\1n\l \ti p(s)\1n\Big\}ds\\
\ns\ds\qq\q -\sum_{j=1}^d\int_t^T\ti q^j(s)dW^j(s),\q t\in[0,T].
\ea\right.\ee
\et

\proof
 The first two equations of \rf{pq- equation} can be unified as follows:
\vskip-3mm
$$\ba{ll}
\ds \3n p(t)\1n=\1nh_x(T)^\top\2n+\1n\dbE_{T-\d}\big[h_y(T)^\top\1n {\bf1}_{[0,T-\d)}(t)\big] \1n+\1n\int_t^T\3n\bigg\{l_x(s)^\top\1n+b_x(s)^\top p(s)\\
\ns\ds \q \1n+\2n\sum_{j=1}^{d}\1n\si_x^j(s)^\top\1n q^j(s)\1n+\1n \ti p(s)\1n+\1n\dbE_s\big[b_y(s\1n+\1n\d)^\top\1n p(s\1n+\1n\d)\1n+\2n\sum_{j=1}^d\si_y^j(s+\d)^\top \1nq^j(s\1n+\1n\d)\1n\\
\ns\ds \q +l_y(s\1n+\1n\d)^\top\2n -\1ne^{-\l\d}\ti p(s\1n+\1n\d)\big]{\bf1}_{[0,T-\d)}(s)\1n\bigg\}ds \1n-\2n\sum_{j=1}^d\1n\int_t^T\3nq^j(s)dW^j(s),  t\in[0,T].
\ea$$
\vskip-1mm
\no For simplicity, in the following, without loss of generality, let $d=1$. By \rf{decompose} and taking the conditional expectation on both sides of \rf{first-order adjoint equation}, it follows that for $0\les t\les T$,
\vskip-6mm
$$\ba{ll}
\ns\ds \dbE_t\big[Y^0(t)+Y^1(t+\d){\bf1}_{[0,T-\d)}(t)\big]=b_x(t)^\top p(t)+\si_x(t)^\top q(t)+l_x(t)^\top+\ti{p}(t)\\
\ns\ds +\dbE_t\big[b_y(t\1n+\1n\d)^\top p(t\1n+\1n\d)+\si_y(t\1n+\1n\d)^\top q(t\1n+\1n\d)+l_y(t\1n+\1n\d)^\top-e^{-\l\d}\ti{p}(t\1n+\1n\d)\big] {\bf1}_{[0,T-\d)}(t),
\ea$$
\vskip-2mm
\no and
\vskip-5mm
$$\ba{ll}
Y^2(t)=b_z(t)^\top p(t)+\si_z(t)^\top q(t)+l_z(t)^\top-\l \ti{p}(t).
\ea$$
\vskip-2mm
\no Noting
\vskip-6mm
$$\ba{ll}
\ns\ds \int_\t^T\3n\dbE_t\bigg[\1n\int_\t^{t+\d}\3n\3nZ^1(t+\d,s)dW(s){\bf1}_{[0,T-\d)}(t)\bigg]dt \1n=\3n\int_\t^{T-\d}\3n\1n\int_\t^t\3nZ^1(t+\d,s)dW(s)dt{\bf1}_{[0,T-\d)}(\t),
\ea$$
\no from \rf{first-order adjoint equation} (c), one has
\vskip-6mm
$$\ba{ll}
\ns\ds \int_\t^T\dbE_t\big[Y^0(t)+Y^1(t+\d){\bf1}_{[0,T-\d)}(t)\big]dt =\int_\t^T\dbE_t\bigg[\dbE_\t\big[Y^0(t)\big]+\int_\t^tZ^0(t,s)dW(s)\\
\ns\ds \qq +\dbE_\t\big[Y^1(t\1n+\1n\d){\bf1}_{[0,T-\d)}(t)\big] +\int_\t^{t+\d}Z^1(t+\d,s)dW(s){\bf1}_{[0,T-\d)}(t)\bigg]dt\\
\ns\ds =\int_\t^T\dbE_\t\big[Y^0(t)+Y^1(t+\d){\bf1}_{[0,T-\d)}(t)\big]dt\\
\ns\ds \q +\int_\t^T\bigg[\int_s^TZ^0(t,s)dt +\int_s^{T-\d}Z^1(t+\d,s)dt{\bf1}_{[0,T-\d)}(s) {\bf1}_{[0,T-\d)}(\t)\bigg]dW(s),
\ea$$
\vskip-1mm
\no and
\vskip-5mm
$$\ba{ll}
\ns\ds \int_\t^TY^2(t)dt=\int_\t^T\bigg[\dbE_\t\big[Y^2(t)\big]+\int_\t^tZ^2(t,s)dW(s)\bigg]dt\\
\ns\ds \qq\qq\q\ =\dbE_\t\bigg[\int_\t^TY^2(t)dt\bigg] +\int_\t^T\int_t^TZ^2(s,t)dsdW(t).
\ea$$
\vskip-2mm
\no Recalling \rf{first-order adjoint equation} (a), one can get
\vskip-6mm
$$\ba{ll}
\ns\ds \eta^0(\t)+\int_\t^T\dbE_\t\[Y^0(t)+Y^1(t+\d) {\bf1}_{[0,T-\d)}(t){\bf1}_{[0,T-\d)}(\t)\]dt +\eta^1(\t){\bf1}_{[0,T-\d)}(\t)\\
\ns\ds +\int_\t^T\bigg\{\z^0(s) +\int_s^T\Big[Z^0(t,s) +Z^1(t+\d,s){\bf1}_{[0,T-\d)}(s) {\bf1}_{[0,T-\d)}(t){\bf1}_{[0,T-\d)}(\t)\Big]dt\\
\ns\ds \qq\q +\z^1(s){\bf1}_{[0,T-\d)}(s)\bigg\}dW(s)\\
\ns\ds =h_x(T)^\top+\int_\t^T\bigg\{b_x(t)^\top p(t)+\si_x(t)^\top q(t)+l_x(t)^\top+\ti{p}(t)+\dbE_t\[b_y(t+\d)^\top p(t+\d)\\
\ns\ds \q +\si_y(t+\d)^\top q(t+\d)+l_y(t+\d)^\top-e^{-\l\d}\ti{p}(t+\d)\] {\bf1}_{[0,T-\d)}(t)\bigg\}dt\\
\ns\ds \q +\dbE_{T-\d}\big[h_y(T)^\top{\bf1}_{[0,T-\d)}(\t)\big],
\ea$$
\vskip-2mm
\no and
$$\ba{ll}
\ns\ds \eta^2(\t)+\int_\t^T\dbE_\t\big[Y^2(t)\big]dt+\int_\t^T\bigg[\z^2(s)+\int_s^TZ^2(t,s)dt\bigg]dW(s)\\
\ns\ds =h_z(T)^\top+\int_\t^T\Big[b_z(t)^\top p(t)+\si_z(t)^\top q(t)+l_z(t)^\top-\l \ti{p}(t)\Big]dt.
\ea$$
\vskip-1mm
\no Thus, the proof is completed.
\endpf

\br
Furthermore, we can arrange the set of anticipated BSDEs \rf{pq- equation} into an anticipated backward SVIE. Since the third equation of \rf{pq- equation} is a linear BSDE, $\ti p(\cd)$ can be expressed by $(p(\cd),q(\cd))$, thus, $(p(\cd), q(\cd))$ satisfies the following anticipated backward SVIE:
\vskip-4mm
\bel{ABSVIEs}\ba{ll}
\ds \3n p(t)\1n=\1nh_x(T)^\top\1n+\dbE_{T-\d} \big[h_y(T)^\top\1n{\bf1}_{[0,T-\d)}(t)\big]\1n +\1n\frac 1 \l \(1\1n-\1ne^{\l((-\d)\vee(t-T))}\1n\)h_z(T)^\top\\
\ns\ds \qq \1n+\3n\int_t^T\3n\bigg\{\1nb_x(s)^\top\1n p(s)\1n+\2n\sum_{j=1}^d\si^j_x(s)^\top\1n q^j(s)\1n+\1nl_x(s)^\top\2n +\1n\dbE_s\bigg[b_y(s\1n+\1n\d)^\top\1n p(s\1n+\1n\d)\\
\ns\ds \qq +\sum_{j=1}^d\si_y^j(s+\d)^\top q^j(s+\d)+l_y(s+\d)^\top\1n\bigg]{\bf1}_{[0,T-\d)}(s)\\
\ns\ds \qq +\dbE_s\bigg[\int_0^{\d\wedge(T-s)}\3n\3n\3ne^{-\l\th} \bigg(b_z(s+\th)^\top\1n p(s+\th)\1n+\2n\sum_{j=1}^d\si_z^j(s\1n+\1n\th)^\top\1n q^j(s\1n+\1n\th) \\
\ns\ds \qq +l_z(s+\th)^\top\bigg)d\th\bigg]\bigg\}ds -\sum_{j=1}^d\int_t^T q^j(s)dW^j(s),\q t\in[0,T].
\ea\ee
\vskip-2mm
\no Now the first-order adjoint equation \rf{first-order adjoint equation} can be unified as the anticipated backward SVIE \rf{ABSVIEs}, and is dual with the variational equations \rf{variational equation-1}-\rf{variational equation-2}, even if the pointwise state delay appears in the terminal cost. This is a new finding.
\er
%
%\proof
%%
%\vskip-2mm
%%
%Substituting \rf{p1q1=pq} into \rf{pq- equation}, then the second equation in \rf{pq- equation} becomes a linear BSDE, thus we have
%%
%\vskip-8mm
%%
%$$\ba{ll}
%%
%\ns\ds \ti p(t)\1n=\1n\dbE_t\1n\int_t^T\3ne^{-\l(s-t)}\bigg(b_z(s)^\top \1n p(s)\1n+\2n\sum_{j=1}^d\si_z^j(s)^\top \1n q^j(s)+\1nl_z(s)^\top\2n\bigg)ds \1n+\1n\dbE_t\[e^{-\l(T-t)}h_z(T)^\top\1n\].
%%
%\ea$$
%%
%\vskip-2mm
%%
%\no It follows that
%%
%\vskip-6mm
%%
%$$\ba{ll}
%%
%\ns\ds \dbE_t\1n\bigg[\1n\int_t^T\ti p(s)ds\1n-\1ne^{-\l\d}\2n\int_t^T\3n\ti p(s+\d){\bf1}_{[0,T-\d)}(s)ds\bigg]\1n=\1n\dbE_t\bigg[\frac 1 \l h_z(T)^\top\1n\(1\1n-\1ne^{\l((-\d)\vee(t-T))}\)\\
%%
%\ns\ds \1n+\2n\int_t^T\3n\int_0^{\d\wedge(T-s)}\3n\3ne^{-\l\th} \bigg(\1nb_z(s+\th)^\top\1n p(s+\th) \1n+\2n\sum_{j=1}^d\si_z^j(s+\th)^\top\1n q^j(s+\th)\1n+\1nl_z(s+\th)^\top\1n\bigg)d\th ds\bigg].
%%
%\ea$$
%%
%\vskip-2mm
%%
%\no Taking the conditional expectation on both sides in \rf{pq- equation} we derive \rf{ABSDE=ABSVIE}.
%%
%\endpf

\br
Let $h_y\equiv0$. Then, by Theorem \ref{theorem relation-1}, $(p(\cd),q(\cd))$ is the unique solution to the following set of anticipated BSDEs:
\vskip-6mm
\bel{pq- ABSDE}\left\{\ba{ll}
\ds \3n-dp(t)=\bigg\{b_x(t)^\top p(t)+\sum_{j=1}^{d}\si_x^j(t)^\top q^j(t)+l_x(t)^\top\2n+\ti{p}(t)+\dbE_t\bigg[l_y(t+\d)^\top\2n \\
\ns\ds \qq\qq+b_y(t+\d)^\top p(t+\d)+\sum_{j=1}^d\si_y^j(t+\d)^\top q^j(t+\d)\1n \\
\ns\ds \qq\qq-e^{-\l\d}\ti{p}(t+\d)\bigg]{\bf1}_{[0,T-\d)}(t)\bigg\}dt-\sum_{j=1}^dq^j(t)dW^j(t),\\
\ds \3n-d\ti{p}(t)\1n=\1n\Big\{b_z(t)^\top\1n p(t)\1n+\2n\sum_{j=1}^d\si_z^j(t)^\top\1n q^j(t)\1n+\1nl_z(t)^\top\2n-\2n\l \ti{p}(t)\Big\}dt\1n-\2n\sum_{j=1}^d\ti{q}^j(t)dW^j(t),\\
\ds \3n p(T)=h_x(T)^\top,\q \ti{p}(T)=h_z(T)^\top.
\ea\right.\ee
\vskip-2mm
%
%\no that satisfies $(p(\cd),q(\cd),\ti{p}(\cd),\ti{q}(\cd))\in L_{\dbF}^2(\Om;C([0,T];\dbR^n))$ $\times L_\dbF^2(0,T;\dbR^{n\times d})$ $\times L_{\dbF}^2(\Om;\\C([0,T];\dbR^n))$ $\times L_\dbF^2(0,T;\dbR^{n\times d})$.
%%
\no Noting \cite{Zhang-2022} assumed that the control domain is convex, then studied the sufficient maximum principle for stochastic optimal control problems with general delays. Let the noisy memory process there disappears, i.e. $X_2^u(\cd)\equiv0$. Then, (10)-(11) in \cite{Zhang-2022} are the same as \rf{pq- ABSDE} above.
\er

%\proof
%%
%\vskip-1mm
%%
%\no By Theorem 4.2 in \cite{Peng-Yang-2009}, ABSDE \rf{pq- ABSDE} has the unique solution. Noting $h_y\equiv0$, we have $(\eta^1(\cd),\z^1(\cd))\equiv(0,0)$, thus it is easy to verify that $(p(\cd),q(\cd),\ti{p}(\cd),\ti{q}(\cd))$ $\in L_{\dbF}^2(\Om;C([0,T];\dbR^n))\times L_\dbF^2(0,T;\dbR^{n})\times L_{\dbF}^2(\Om;C([0,T];\dbR^n))\times L_\dbF^2(0,T;\dbR^{n})$ defined by \rf{pq-}, furthermore by  Theorem \ref{theorem relation-1}, $(p(\cd),q(\cd),\ti{p}(\cd)$, $\ti{q}(\cd))$ defined by \rf{pq-} satisfies \rf{pq- ABSDE}, thus the proof is completed.
%%
%\endpf

\br
Let $h_y,h_z\1n\equiv\1n0$. Then, \rf{ABSVIEs} becomes a simple anticipated BSDE consistent with (5.1) in \cite{Meng-Shi-2021}, when the distributed delay disappears in Problem {\rm (P)}.
%
%\vskip-5mm
%%
%\bel{ABSDE-}\left\{\ba{ll}
%%
%\ds \3n -dp(t)=\bigg\{b_x(t)^\top p(t)+\sum_{j=1}^d\si^j_x(t)^\top q^j(t)+l_x(t)^\top +\dbE_t\bigg[l_y(t+\d)^\top\\
%%
%\ns\ds \qq\qq +b_y(t+\d)^\top p(t+\d)+\sum_{j=1}^d\si_y^j(t+\d)^\top q^j(t+\d)\bigg]{\bf1}_{[0,T-\d)}(t) \\
%%
%\ns\ds \qq\qq
%  +\dbE_t\bigg[\int_0^{\d\wedge(T-t)}\2n e^{-\l\th}\bigg(l_z(t+\th)^\top+b_z(t+\th)^\top p(t+\th)\\
%  %
%  \ns\ds \qq\qq\1n+\2n\sum_{j=1}^d\si_z^j(t+\th)^\top\1n q^j(t+\th)\bigg)d\th\bigg]\bigg\}dt\1n-\2n\sum_{j=1}^dq^j(t)dW^j(t), t\in[0,T],\\
%%
%\ns\ds \3n p(T)=h_x(T)^\top,
%%
%\ea\right.\ee
%%
\er
%
%\proof
%%
%\vskip-1mm
%%
%This theorem is easy to be obtained by $h_z\equiv 0$ and Theorem \ref{theorem relation-2}.
%%
%\endpf

\subsubsection{Extensions of second-order adjoint equations}

In this subsection, we study three typical control systems to demonstrate second-order adjoint equations more clearly.

\ms

\textbf{Case I: Stochastic optimal control problems without delay}

\ss

In this case, Problem (P) becomes a classical stochastic optimal control problem. From \rf{P1}, \rf{P2-1}, \rf{P3-1} and \rf{P4}, $P_1^{(11)}(r)$, $P_2^{(11)}(r)$, $P_3^{(11)}(r)$, $P_4^{(11)}(\th,r)\neq0$, $0\les r,\th\les T$, and other terms in \3n \rf{second-order adjoint equations} are all 0. Then, \rf{cP} becomes
\vskip-5mm
$$\ba{ll}
\ns\ds \cP^1(s)\equiv\cP(s) =\1nh_{xx}(T)+\int_s^T \3n\dbE_r\bigg[b_x(r)^\top\cG_2^{(1)}(r)\1n+\2n\sum_{j=1}^d \si_x^j(r)^\top\cQ_{2j}^{(1)}(r)\1n+\1n\cG_2^{(1)}(r)^\top b_x(r)\\
\ns\ds \q \1n+\2n\sum_{j=1}^d\1n\cQ_{2j}^{(1)}(r)^\top\1n\si_x^j(r)\1n\bigg]\1ndr \2n+\3n\int_s^T\3n\bigg\{\1nb_x(r)^\top\1n\dbE_r\1n\bigg[\2n\int_r^T \3n\2n\cG_4^{(1)}(\th,r)d\th\bigg] \2n+\2n\sum_{j=1}^d\si_x^j(r)^\top\3n\int_r^T\3n \1n \cQ_{4j}^{(1)}(\th,r)d\th\bigg\}dr\\
\ns\ds \q +\int_s^T\bigg\{\bigg(\int_r^T \dbE_{r}\big[\cG_4^{(1)}(\th,r)^\top\big] d\th\bigg)b_x(r) +\sum_{j=1}^d\bigg(\int_r^T \cQ_{4j}^{(1)}(\th,r)^{\top}d\th \bigg) \si_x^j(r)\bigg\}dr\\
\ns\ds \q +\int_s^T \bigg\{G_{xx}(r)+\sum_{j=1}^d\si_{x}^j(r)^\top\dbE_r\big[\cP(r)\big]\si_x^j(r)\bigg\}dr.
\ea$$
\vskip-2mm
\no Denote
\vskip-4mm
\bel{case1-bar P}\ba{ll}
\ns\ds \bar P^1(\t):=\dbE_\t\big[\cP^1(\t)\big],\q \bar Q^1(\t):=\cQ_2^{(1)}(\t)+\int_\t^T\cQ_4^{(1)}(\th',\t)d\th'.
\ea\ee
\vskip-1mm
\no Then, $(\bar P^1(\cd),\bar Q^1(\cd))$ satisfies the following BSDE:
\vskip-4mm
\bel{case1-BSDE}\ba{ll}
\ns\ds \bar P^1(\t)=h_{xx}(T)+\int_\t^T\bigg\{b_x(t)^\top\bar P^1(t)+\sum_{j=1}^d\si_x^j(t)^\top\bar Q_j^1(t)+\bar P^1(t)^\top b_x(t)\\
\ns\ds \q +\sum_{j=1}^d\bar Q_j^1(t)^\top\si_x^j(t)+l_{xx}(t)+\blan p(t),b_{xx}(t)\bran +\sum_{j=1}^d\blan q^j(t),\si_{xx}^j(t)\bran\\
\ns\ds \q +\sum_{j=1}^d\si_x^j(t)^\top\bar P^1(t)\si_x^j(t)\bigg\}dt-\int_\t^T\sum_{j=1}^d\bar Q_j^1(t)dW^j(t),\q \t\in[0,T],\\
\ea\ee
\vskip-2mm
\no which is consistent with (3.9) in \cite{Yong-Zhou-1999}.

In fact, we have
\vskip-5mm
$$\ba{ll}
\ns\ds\cP^1(\t)=\cG_2^{(1)}(\t)+\int_\t^T \cG^{(1)}_4(\th',\t)d\th',
\ea$$
\vskip-2mm
\no and by \rf{case1-bar P},
\vskip-5mm
\bel{case1-} \ba{ll}
\ns\ds \cP^1(\t)=\1nh_{xx}(T)+\int_\t^T\bigg\{b_x(t)^\top\bar P^1(t)+\sum_{j=1}^d\si_x^j(t)^\top\bar Q_j^1(t)+\bar P^1(t)^\top b_x(t)\\
\ns\ds \qq\qq +\sum_{j=1}^d\bar Q_j^1(t)^\top\si_x^j(t)+\blan p(t),b_{xx}(t)\bran+\sum_{j=1}^d\blan q^j(t),\si_{xx}^j(t)\bran \\
\ns\ds \qq\qq+l_{xx}(t)+\sum_{j=1}^d\si_x^j(t)^\top\bar P^1(t)\si_x^j(t)\bigg\}dt.
\ea\ee
\vskip-1mm
\no By the first equation of \rf{second-order-adjoint-constraints}, we have
\vskip-4mm
\bel{case1--}\ba{ll}
\ns\ds h_{xx}(T)=\dbE_\t\big[h_{xx}(T)\big] +\sum_{j=1}^d\int_\t^TQ_{1j}^{(11)}(r)dW^j(r),\q \t\in[0,T].
\ea\ee
\vskip-1mm
\no Noting \rf{second-order-adjoint-constraints}, for $k=2,3,\ i,\ell=1,2,3$, we get
\bel{case1---}\ba{ll}
\ns\ds \int_\t^T\3nP_k^{(i\ell)}(\th)^\top \1nd\th\1n=\1n\dbE_\t \bigg[\1n\int_\t^T\3nP_k^{(i\ell)}(\th)^\top\1n d\th\bigg] \1n+\2n\sum_{j=1}^d\1n\int_\t^T\3n\int_\t^\th\3n Q_{kj}^{(i\ell)}(\th,r)^\top dW^j(r)d\th\\
\ns\ds =\dbE_\t\bigg[\int_\t^TP_k^{(i\ell)}(\th)^\top d\th\bigg] +\sum_{j=1}^d\int_\t^T\int_r^T Q_{kj}^{(i\ell)}(\th,r)^\top d\th dW^j(r),
\ea\ee
\vskip-4.2mm
\no and
\vskip-4mm
\bel{case1----}\ba{ll}
\ns\ds \int_\t^T\int_\t^TP_4^{(i\ell)}(\th',\th)d\th d\th'=\int_\t^T\int_\t^T\dbE_{\t}\big[P_4^{(i\ell)}(\th',\th)\big]d\th d\th'\\
\ns\ds \q +\sum_{j=1}^d\int_\t^T\int_r^T\int_r^TQ_{4j}^{(i\ell)}(\th',\th,r)d\th d\th' dW^j(r),\q i,\ell=1,2,3.
\ea\ee
\vskip-1mm
\no From \rf{case1-}-\rf{case1----}, we obtain
\vskip-5.5mm
$$\ba{ll}
\ns\ds \cP^1(\t)=\bar P^1(\t)+\sum_{j=1}^d\int_\t^T\bar Q^1_j(t)dW^j(t),\q \t\in[0,T],
\ea$$
\vskip-2mm
\no which implies \rf{case1-BSDE}.

\ms

\textbf{Case II: Stochastic optimal control problems with control delay only}

\ms

\no In this case, $b_y,b_z,\si_y,\si_z,l_y,l_z,h_y,h_z\1n=\1n0$. From \rf{P1},\1n \rf{P2-1},\1n \rf{P3-1}\1n and \1n\rf{P4}, we have
\vskip-3mm
$$\ba{ll}
\ns\ds P_1^{(11)}(r),\ P_2^{(11)}(r),\ P_3^{(11)}(r),  P_4^{(11)}(\th,r),\ P_4^{(12)}(\th,r)\neq0,\q 0\les r\les \th\les T,\\
\ns\ds P_4^{(11)}(\th,r),\ P_4^{(21)}(\th,r)\neq0,\q 0\les\th<r\les T,

\ea$$
\vskip-2mm
\no and other terms in \rf{second-order adjoint equations} are all 0. From \rf{psi4-ij} and \rf{g4-ij}, we obtain
\vskip-8mm
$$\ba{ll}
\ns\ds \psi_{4}^{(12)}(\th,r)=0,\q g_{4}^{(12)}(\th,\th',r)=b_x(r)^\top P_4^{(12)}(\th,\th')+\sum_{j=1}^d\si_x^j(r)^\top Q_{4j}^{(12)}(\th,\th',r),
\ea$$
\vskip-3mm
\no and then,
\vskip-4mm
\bel{case2-P3-12-1}\ba{ll}
\ns\ds P_4^{(12)}(\th,r)=\int_r^T\bigg[b_x(r)^\top P_4^{(12)}(\th,\th')+\sum_{j=1}^d\si_x^j(r)^\top Q_{4j}^{(12)}(\th,\th',r)\bigg]d\th'\\
\ns\ds \qq\qq\qq -\sum_{j=1}^d\int_r^T Q_{4j}^{(12)}(\th,r,\th')dW^j(\th'),\q \th\ges r.
\ea\ee
\vskip-1mm
\no On the other hand, recalling \rf{second-order-adjoint-constraints}, for $\th\ges r$ , we have
\vskip-4mm
\bel{case2-P3-12-2}\ba{ll}
\ns\ds P_4^{(12)}(\th,r)=\dbE_{\th'}\big[P_4^{(12)}(\th,r)\big] +\sum_{j=1}^d\int_{\th'}^rQ_{4j}^{(12)}(\th,r,s)dW^j(s).
\ea\ee
\vskip-1mm
\no By the unique solvability of the backward SVIEs, \rf{case2-P3-12-1} and \rf{case2-P3-12-2} lead to that
\vskip-4mm
$$\ba{ll}
\ns\ds P_4^{(12)}(\th,r)=0,\q Q_4^{(12)}(\th,r,\th')=0,\qq \th\ges r.
\ea$$
\vskip-1mm
\no Hence, it follows that for $\th\ges r$,
\vskip-4mm
$$\ba{ll}
\ns\ds \cG_4^{(2)}(\th,r)=\int_r^TP_4^{(12)}(\th,\th')d\th'=0,\\
\ns\ds \cQ_4^2(\th,r)\1n=\3n\int_r^T\3n\1nQ_4^{(12)}(\th,\th',r)d\th' \1n=\2n\int_r^\th \3n Q_4^{(12)}(\th,\th',r)d\th'\1n+\2n\int_\th^T\3nQ_4^{(21)} (\th',\th,r)^\top\1n d\th'\1n=\1n0.
\ea$$
\vskip-3mm
\no Then, \rf{cP} becomes
\vskip-4mm
$$\ba{ll}
\ns\ds \cP^2(s)\equiv\cP(s) =h_{xx}(T)+\int_s^T \dbE_r\bigg[b_x(r)^\top\cG_{2}^{(1)}(r) +\sum_{j=1}^d \si_x^j(r)^\top\cQ_{2j}^{(1)}(r)\\
\ea$$
\vskip-7mm
$$\ba{ll}
\ns\ds \q +\cG_2^{(1)}(r)^\top b_x(r) +\sum_{j=1}^d\cQ_{2j}^{(1)}(r)^\top\si_x^j(r)\bigg]dr +\int_s^T\bigg\{b_x(r)^\top\dbE_r\bigg[\int_r^T \cG_4^{(1)}(\th,r)d\th\bigg]\\
\ns\ds \q +\sum_{j=1}^d\si_x^j(r)^\top\int_r^T \cQ_{4j}^{(1)}(\th,r)d\th\bigg\}dr +\int_s^T\bigg\{\bigg(\int_r^T\dbE_{r}\big[\cG_4^{(1)}(\th,r)^\top\big]d\th\bigg)b_x(r)\\
\ns\ds \q +\1n\sum_{j=1}^d\bigg(\1n\int_r^T\3n \cQ_{4j}^{(1)}(\th,r)^\top\1n d\th\bigg)\si_x^j(r)\1n\bigg\}dr \1n+\1n\int_s^T \3n \bigg\{G_{xx}(r)\1n+\2n\sum_{j=1}^d\si_{x}^j(r)^\top\1n \dbE_r\big[\cP(r)\big]\si_x^j(r)\1n\bigg\}dr.
\ea$$
\vskip-2mm
\no Denote
\vskip-4mm
$$\ba{ll}
\ns\ds \bar P^2(\t):=\dbE_\t\big[\cP^2(\t)\big],\q \bar Q^2(\t):=\cQ_2^{(1)}(\t)+\int_\t^T\cQ_4^{(1)}(\th',\t)d\th'.
\ea$$
\vskip-1mm
\no Then, similar to Case I, $(\bar P^2(\cd),\bar Q^2(\cd))$ also satisfies the BSDE \rf{case1-BSDE}.

\ms

\textbf{Case III: Linear quadratic stochastic optimal control problems}

\ss

Consider the following state equation:
\vskip-4mm
$$\left\{\ba{ll}
\ds dX(t)=\big[A(t)X(t) +B(t)u(t)+\bar{B}(t)u(t-\delta)\big]dt\\
\ns\ds\qq\qq +\big[\bar{C}(t)X(t-\delta) +D(t)u(t)+\bar{D}(t)u(t-\delta)\big]dW(t),\ t\in[0,T],\\
\ns\ds X(t)=\xi(t),\ u(t)=\eta(t),\ t\in[-\delta,0],
\ea\right.$$
\vskip-2mm
\no with the quadratic cost functional
\vskip-7mm
$$\ba{ll}
\ns\ds J(u(\cdot))=\dbE\big[\blan GX(T),X(T)\bran+2\blan g,X(T)\bran\big]\\
\ns\ds \q +\mathbb{E}\2n\int_0^T\3n\left\langle
  \left[\3n\begin{array}{cccc}
  Q_{00}(t) & 0 & S_{00}(t)^\top & S_{01}(t)^\top\\
  0 & Q_{11}(t) & S_{10}(t)^\top & S_{11}(t)^\top \\
  S_{00}(t) & S_{10}(t) & R_{00}(t) & R_{01}(t)\\
  S_{01}(t) & S_{11}(t) & R_{01}(t)^\top & R_{11}(t)
  \end{array}\3n\right]
  \left[\3n\begin{array}{c}X(t) \\ X(t-\d) \\ u(t)\\ u(t-\d)\end{array}\3n\right],
  \left[\3n\begin{array}{c}X(t) \\ X(t-\d) \\ u(t)\\ u(t-\d)\end{array}\3n\right]
  \right\rangle dt,
\ea$$
\vskip-2mm
\no where
$A(\cdot),B(\cdot),\bar{B}(\cdot),\bar{C}(\cdot),D(\cdot), \bar{D}(\cdot),Q_{00}(\cdot)$, $S_{00}(\cd),S_{01}(\cd)$, $Q_{11}(\cd),S_{10}(\cd),S_{11}(\cd)$, $R_{00}(\cdot),R_{01}(\cdot)$, $R_{11}(\cd)$ are all deterministic functions,
and $G\in\dbR^{n\times n}$, $g\in\dbR^n$. In this case, \rf{P1}, \rf{P2-1}, \rf{P3-1} and \rf{P4} become
\vskip-7mm
$$\ba{ll}
\ns\ds P_1^{(11)}(r)=G,\q P_2^{(11)}(r)=A(r)^\top\bigg[P_1^{(11)}(r)+\int_r^TP_2^{(11)}(\th)d\th\bigg],\q 0\les r\les T,\\
\ns\ds P_4^{(11)}(\th,r)=A(r)^\top\cG_4^{(1)}(\th,r),\q P_4^{12}(\th,r)=A(r)^\top\cG_4^{(2)}(\th,r),\q 0\les r\les\th\les T,\\
\ns\ds P_4^{(11)}(\th,r)=\cG_4^{(1)}(r,\th)^\top A(\th),\q P_4^{(21)}(\th,r)=\cG_4^{(1)}(r,\th)^\top A(\th),\q 0\les\th< r\les T,\\
\ns\ds P_3^{(11)}(r)=Q_{00}(r),\q P_3^{(22)}(r)=Q_{11}(r)+\bar C(r)^\top\cP(r)\bar C(r),\q 0\les r\les T,
\ea$$
\vskip-2mm
\no and other terms in \rf{second-order adjoint equations} are all 0. From \rf{cG3-CQ3} we have
\vskip-7mm
$$\ba{ll}
\ns\ds \cG_4^{(1)}(\th,r)\1n=\1nP_2^{(11)}(\th)^\top\1n+P_3^{(11)}(\th) \1n+\1n\int_r^T\1n\[P_4^{(11)}(\th,\th') \1n+\1n{\bf1}_{(\d,\i)}(\th'-r)P_4^{(21)}(\th,\th')\]d\th',\\
\ns\ds \cG_4^{(2)}(\th,r)={\bf1}_{(\d,\i)}(\th-r)P_3^{(22)}(\th) +\int_r^TP_4^{(12)}(\th,\th')d\th'.
\ea$$
\vskip-2mm
\no Let $\th-\d\les r\les\th$, $\t+\d\les\th\les T$, and consider
\vskip-7mm
$$\ba{ll}
\ns\ds \cG_4^{(2)}(\th,r)=\int_r^TP_4^{(12)}(\th,\th')d\th' =\int_r^{\th}P_4^{(12)}(\th,\th')d\th'=\int_r^\th A(\th')^\top\cG_4^{(2)}(\th,\th')d\th'.
\ea$$
\vskip-2mm
\no Then, we have
\vskip-4mm
$$\ba{ll}
\ns\ds \cG_4^{(2)}(\th,r)=0,\q  \th-\d\les r\les\th,\q \t+\d\les\th\les T.
\ea$$
\vskip-2mm
\no Hence, \rf{cP} becomes
\vskip-5.5mm
$$\ba{ll}
\ds \cP^3(s)\equiv\cP(s)=G+\int_s^T\big[A(r)^\top\cG_2^{(1)}(r)+\cG_2^{(1)}(r)^\top A(r)\big]dr\\
\ns\ds \q +\int_s^TA(r)^\top\Big(\int_r^T\cG_4^{(1)}(\th,r)d\th+\int_{r+\d}^T\cG_4^{(2)}(\th,r)d \th\Big)dr\\
\ns\ds \q +\int_s^T\Big(\int_r^T\cG_4^{(1)}(\th,r)^\top d\th +\int_{r+\d}^T\cG_4^{(2)}(\th,r)^\top d\th\Big)A(r)dr +\int_s^TQ_{00}(r)dr\\
\ns\ds \q +\int_{s+\d}^T \big[Q_{11}(r)+\bar C(r)^{\top}\cP(r)\bar C(r)\big]dr.
\ea$$
\vskip-1.5mm
\no Similar to Case I, $\cP^3(\cd)$ satisfies the following ordinary differential equation:
\vskip-5mm
\bel{case3-ODE}\left\{\ba{ll}
\ds -\dot\cP^3(s)=A(s)^\top\cP^3(s)+\cP^3(s)^\top A(s)+Q_{00}(s)+\big[Q_{11}(s+\d)\\
\ns\ds \qq\qq\ +\bar C(s+\d)^\top\cP^3(s+\d)\bar C(s+\d)\big]{\bf1}_{[0,T-\d)}(s),\q \ae s\in[0,T],\\
\ns\ds \cP^3(T)=G.
\ea\right.\ee

\br
\vskip-1.5mm
For Case I, when the delay disappears in the control system, the equation \rf{ABSVIEs} satisfied by $(p(\cd),q(\cd))$, becomes (3.8) in \cite{Yong-Zhou-1999}; the equation \rf{case1-BSDE} satisfied by $\cP(\cd)$, becomes (3.9) in \cite{Yong-Zhou-1999}, and so Theorem \ref{MP} reduces to Theorem 3.2 in \cite{Yong-Zhou-1999}. For Case II and Case III, \rf{ABSVIEs}, \rf{case1-BSDE} and \rf{case3-ODE} are consistent with (5.1) and (5.2) in \cite{Meng-Shi-2021}, respectively, thus Theorem \ref{MP} reduces to Theorem 5.1 in \cite{Meng-Shi-2021}.
\er

\section{Concluding remarks}

In this paper, a stochastic optimal control problem is considered and the control domain is allowed to be non-convex. The pointwise state delay, distributed state delay and pointwise control delay can appear in the diffusion term and the terminal cost. Via the theory of backward stochastic Volterra integral systems, we transform delayed variational equations into Volterra integral equations without delay, introduce some new second-order adjoint equations and derive a general maximum principle, without any additional conditions. Finally, to express adjoint equations more compact, we in detail discuss them for three typical control systems.

\end{document}